
\font\steptwo=cmb10 scaled\magstep2
\font\stepthree=cmb10 scaled\magstep4
\magnification=\magstep1
\settabs 18 \columns
 
\hsize=16truecm
\baselineskip=17 pt

\def\b{\bigskip}
\def\bb{\bigskip\bigskip}

\def\no{\noindent}
\def\r{\rightline}
\def\ce{\centerline}
\def\ve{\vfill\eject}

\def\r{\rightline}
\font\got=eufm8 scaled\magstep1

\def\~{\hskip-1mm}
\def\g{{\got g}}

\def\harr#1#2{\smash{\mathop{\hbox to .25 in{\rightarrowfill}}
 \limits^{\scriptstyle#1}_{\scriptstyle#2}}}

\font\aab=cmbsy10 at 18pt
\def\bigbullet{\lower.3ex\hbox{{\aab  \char'017}}}

\def\today{\ifcase\month\or January\or February\or March\or
April\or  May\or June\or July\or August\or September\or
October\or November\or  December\fi \space\number\day,
\number\year }

\r \today

\b

\def\Rit{\hbox{\it I\hskip -2pt R}}

\def\Cit{\hbox{\it l\hskip -5.5pt C\/}}
\def\Rrm{\hbox{\rm I\hskip -2pt R}}
\def\Nrm{\hbox{\rm I\hskip -2pt N}}
\def\Crm{\hskip0.5mm \hbox{\rm l\hskip -5.5pt C\/}}

 \bb\bb  {\ce {\stepthree    Deformation Quantization}{\ce {\stepthree 
on Singular Coadjoint Orbits.}}}
\b

 {\ce {Christian Fr\o nsdal }}

 \def\g{\hbox{\got g}}

 {\it\ce {Physics Department, University of California,} \ce{Los Angeles
CA 90024, USA, fronsdal@physics.ucla.edu}} 
\b
\no{\it  ABSTRACT.}~ Invariant star products are constructed on
minimal coadjoint orbits of all the simple Lie algebras.  Explicit expressions
are given for the generators of the Joseph ideals and the associated
infinitesimal characters.
\b
   
\no{ \bf Content:}~

{\bf 1. Introduction.} History, cohomology, results, outline.

{\bf 2. Associative *-products   
and cohomology.} ~~~~~ 2.1. Formal star-products. 

 \hskip.5cm 2.2. The BGS 
decomposition of 
 Hochschild (co-)homology.   2.3. Star products

  \hskip.5cmand the  BGS
decomposition.

{\bf 3. Some varieties with singularities.} 3.1. Conic varieties defined by
quadratic

 \hskip.5cm relations. 
~~~ 3.2. The simple cone; one  quadratic relation. ~~~~3.3.  Case of one

 \hskip.5cmpolynomial relation.

{\bf  4. Invariant star products on coadjoint orbits.} ~~ 4.1.
Background.

 \hskip.5cm 4.2. Invariant star products.~ ~4.3. Finite dimensional
representations.

{\bf  5. Introduction to singular orbits.}~~~ 5.1. Coadjoint orbits of
so(2,1). 

 \hskip.5cm 5.2.   Minimal orbits of 
 sl($n$).    
 ~5.3. Minimal orbits of sp(2n), so(n) and   

 \hskip.5cm  the others.  
 ~~5.4. Associated representations and Joseph ideals. 

{\bf 6. Invariant star products on minimal orbits.}
 6.1. Computation of the 

 \hskip.5cmhomology. 
~~6.2. Invariant star products and the correspondence principle. 

 \hskip.5cm6.3.
Calculations for sl($n$); Joseph ideal and 
highest weight modules. 

 \hskip.5cm6.4 Calculations for so($n$). 6.5 Highest
weight module for so(2l+1) and for so(2l). 

 \hskip.5cm 6.6.
Uniform calculations for the exceptional, simple Lie algebras. ~~Appendix.

\ve  
 
\no{\steptwo   1. Introduction.}

Let \g~ be a Lie algebra over a field $K$, \g$'$ its K-vector space
dual and
$G$ its adjoint Lie group. The group acts smoothly on \g$'$ and each orbit is a
variety of the type $K[x_1,...,x_N]/(R)$, where
$R$ is a set of polynomial relations.  And each variety in this very rich
collection  comes fitted with a natural symplectic structure. As shown by
Kirillov, there are interesting relationships between these coadjoint orbits
and representation theory; Kostant [K1] and Souriau [S] brought in ideas from
classical and quantum mechanics and especially from quantization. Meanwhile
Gerstenhaber had developed his deformation theory [G1], and eventually it was
understood that quantization on a symplectic space is a deformation of the
algebra of functions. One can go further, to regard every deformation as  
quantization.  

It is remarkable, we think, that the greater part of recent
work on this subject has tended to downplay the special features that
coadjoint orbits inherit from the Lie algebra. This paper is a contribution to
the study of invariant star products, a type of deformation-quantization that
makes more intimate contact with Lie structure. We shall also borrow 
from the analogy with quantization of
mechanical and field theoretic systems. It will turn out that the view of
quantization that was introduced by Weyl, and that is intimately tied to
Bohr's correspondence principle, is a most efficient one; it supplements
cohomological methods precisely where those become inefficient. 

 \b
\ce{\bf  History} 

A first correspondence between quantum mechanical and classical
observables was established by Weyl [Wl] (classical-to-quantum) and
Wigner [Wr] (the other way). The formula for pulling back the
non-commutative commutator of operators to classical observables is due to
Moyal [M].   The resulting
deformation of the Poisson bracket, and of the product of functions on
phase space, was subjected to mathematical analysis by Vey [Vy], following
the study of one-differentiable deformations in [FLS].

The papers [BFFLS] generalized this notion and proposed a new axiomatic
approach to quantization,  interpreted as a deformation  and formulated  in
terms of a general type of associative star product. Following Moyal and
Vey one considers the space 
 of formal series, in a parameter $\hbar$, of $C^\infty$ functions on a
symplectic space $W$, with an associative product of the form
$$
f*g = fg + {i\hbar\over 2}\{f,g\} + \sum _{n = 2}^\infty \hbar^nC_n(f,g),~~
f,g
\in C^\infty W,
$$
where $fg$ is the ordinary product of functions, $\{.,.\}$ is the
Poisson bracket, and the cochains $C_n$ are often (but not in
this paper) taken to be bi-differential operators.  There were at least four
different developments.

1. The problem of classification of
star products on symplectic spaces, up to a natural but very weak form of
equivalence, was investigated by Gutt and others  [Gt]. The existence of
star products on an arbitrary symplectic manifold was established by
de Wilde and Lecomte [dL] and by Fedosov  [F], culminating with the results
of Kontsevich  who demonstrated the existence of star products on an arbitrary
Poisson manifold  [Kh],[T].  The importance of these results is that they are
global statements about (smooth) manifolds.  

2. A  generalization of Weyl's original correspondence is required in
field theories and is referred to as `the ordering  problem'. There have been
applications to mechanical problems as well [AW]. Attention is called
to the  power of this method in the local algebraic context. From this
point of view {\it the existence} of associative star products is  not
surprising, although a concise expression for $f*g$
may be difficult to obtain.

3. The concept of invariant star products on coadjoint orbits of Lie
algebras was proposed in [BFFLS]. The existence of invariant star products on
any `regular'  coadjoint  orbit of a semisimple Lie algebra was demonstrated in
the same paper. This result was obtained by setting up an explicit, invariant
type of generalized Weyl correspondence. The original Weyl correspondence
yields, in particular, an invariant star product  for the linear symplectic
algebra of the manifold. Invariant star products are implicit in recent
studies of nilpotent orbits, especially those that deal with the Joseph ideal,
e.g. [BJ]. There was an important parallell development in the work of Berezin
[Bn].

4. Star quantization was used as a tool in representation
theory, to generalize the method of geometric quantization of Kostant [K1]
and Souriau [S]. See for example [F1], [Gt], [W]. This idea
has not yet realized its full potential. 
 \ve
\ce{\bf Cohomology}

Since the first paper by Vey it has been clear that the classification of
deformations of the algebra of functions of a smooth manifold is intimately
related to the Hochschild cohomology of the manifold. Less well known is the
role played by the BGS decomposition of this complex. As an example we recall
that the nonexistence of abelian deformations on a smooth manifold is related
to the vanishing of the Harrison component of Hochchild cohomology. More
general algebraic manifolds offer more room for deformations.   
 The lifting of a Poisson structure to a star product is governed by
components of cohomology that are purely local, being associated with the
singularities. This strongly suggests that  global Poisson structures on  
varieties more general than smooth manifolds lift to global star products.
But the results reported here are local.
\b

\ce{\bf  Results}
This paper was intended as a preliminary study of the deformations of the
coordinate algebras of some algebraic varieties with singularities, a context
in which the BGS decomposition could be expected to have some
interesting applications.  The coadjoint orbits of simple Lie groups offers an 
especially rich and interesting family of examples. 

A familiar reduction paradigm was used to reduce the cohomology to a complex
of closed, linear chains. It had been expected that this would lead to an
easy classification of essential deformations.
It  turned out, however, to be difficult to obtain enough information about
the homology of the reduced complex. Before us experts [BJ], applying much
heavier machinery ([Be], [BGS]) to the same problem, had the same experience
and were forced to fall back on  heuristic (?) arguments.  See below, and in
Section 6.3.   

Returning to the point of view that sees a star product as a correspondence
between ordinary polynomials and star polynomials, we were able to 
complete the calculations. The principal results of this paper
are as follows.

  (a) The  Hochschild cohomology of the coordinate
algebra  of an algebraic variety defined by quadratic relations   is
isomorphic to that of its restriction to linear, closed chains. Obstructions
to extending a first order star product to a formal, associative product, to
all orders in the deformation parameter, can therefore be reduced to a study of
the star products  $x*x$ and
$x*x*x$ for $x$ of degree 1. 
\ve

 (b) The minimal, coadjoint orbit of a complex, simple Lie algebra $\neq$
sl($n$) admits a one-parameter  
 family of invariant star    products/deformations. In the case of sl($n$)
we determine a 2-parameter family of deformations, including an interesting 
1-parameter abelian subfamily.  In all the other cases there is a unique,
invariant star product such that,
$\forall u,v
\in
\g\,'$,
$$
u*v -v*u = \hbar\{u,v\},~~ u*u = uu + k,
$$
for some $k\in K$; see Section 6. The value of $k$ is determined by an
examination of the next case,
$$
u*u*u = uuu + \phi(u),
$$
where $\phi$ is a polynomial of degree one. Both $k$ and $\phi$ are
uniquely determined by the relations that define the orbit.  A
uniform calculation covers the five exceptional algebras. 
For sl($n$) and the other classical simple Lie algebras we calculate the
generators of the Joseph ideals and determine the associated highest weight
modules. The uniqueness of the Joseph ideal is a corollary, see [WS].
 \b
\ce{ Outline.}
Section 2 contains a short introduction to formal star products and the BGS
decomposition of Hochshild cohomology. Section 3 is a study of the Hochschild
cohomology of  algebraic varieties defined by a set of quadratic relations.
  A principal tool used here is a
reduction of the Hochschild complex to a subcomplex of closed, linear chains;
Theorem 3.1.2. The method is effective when the underlying algebra is
finitely generated and thus graded, with only positive degrees. The BGS
decomposition is  used throughout. 

Section 4 is a brief introduction to invariant star products. Section 4.3 gives
an example of the appearance of finite representations within the program of
star quantization. Included here is the first example (known to me) of harmonic
polynomials in the enveloping algebra of a simple Lie algebra. We show an
`instance of a deformation' (not a formal deformation) in which the deformed
variety (the spectrum of the deformation of the ring of coordinate functions)
is a finite union of disconnected varieties. 

Section 5 is an informal study of singular, nilpotent orbits. It serves to
introduce this subject to nonexperts, and to build some support for our own
intuition. Associated to these orbits, and to the Joseph ideals, are certain 
very special, unitary representations that play a conspicuous role in physics.

Section 6 examines invariant star products on the most interesting coadjoint
orbits, those of minimal dimensions, with their  Joseph ideals.
Attempting to calculate the cohomology we encounter a difficulty that had
already been met by Braverman and Joseph  [BJ],  and fail to obtain a
sufficiently detailed  description of the space Hoch$_3$ of the coordinate
algebra. For the solution of this problem we offer only 
conjecture  6.1.1, but we circumvent  the difficulty by an independent, direct
calculation. It is done by reformulating the search for an invariant
star product as a correspondence principle, in the spirit of Weyl's symmetric
ordering. Detailed knowledge of the cohomology of the restricted complex is
not needed. Generators of the Joseph ideals are determined.

The Lie algebras sl($n$) and so($n$) are treated separately and all the
calculations are included, with proofs relegated to an Appendix. The case of
sp(2$n$) is too well known to warrant much attention. The five exceptional
simple Lie algebras are handled uniformly together, all the calculations are
in the main text,  they are both short and easy.

Within the family of generally noncommutative star products there may be a
subfamily of non trivial abelian ones.  According to Braverman and Joseph, this
would contradict  the fact that the minimal orbits - excepting the case of  
sl(n) - are rigid. Granted that a deformation of the ring of coordinate
functions implies a deformation of its spectrum; but is it known that the
deformed spectrum of an abelian deformation is always an algebraic variety, or
that every equivariant deformation of a coadjoint orbit is a coadjoint orbit?
In any case we confirm that abelian deformations of the coordinate algebra
exists only in the case of sl($n$), and that in that case the spectrum can be
identified with a neighbouring orbit of the same dimension.

\b\b\b

\b \b
\no{\steptwo  2. Associative  star products and cohomology.}~
 
\ce{\bf  2.1 Formal   $*$-products.}~ A formal, abelian
$*$-product on a commutative algebra $A$ ~  is a commutative,
associative product on the space of formal power series in a formal
parameter
$\hbar$  with coefficients in $A$, given by a formal series
$$
f*g = fg + \sum_{n>0} \hbar^nC_n(f,g).\eqno(2.1)
$$
Associativity is the condition that $f*(g*h) = (f*g)*h $, or
$$
\sum_{m,n = 0}^k \hbar^{m+n}\biggl(
C_m(f,C_n(g,h))-C_m(C_n(f,g),h))\biggr) = 0,\eqno(2.2)
$$
where $C_0(f,g) = fg$. This must be interpreted as an identity in
$\hbar$; thus
$$
\sum_{m,n = 0}^k \delta_{m+n,k}\biggl(
C_m(f,C_n(g,h))-C_m(C_n(f,g),h) )\biggr) = 0, ~~k =
1,2,\cdots~.\eqno(2.3)
$$
The  formal $*$-product (2.1) is associative to order $p$ if Eq.(2.3)
holds for $k = 1,\cdots p$.

A first order   $*$-product is a product
$$
f*g = fg + \hbar C_1(f,g), \eqno(2.4)
$$
associative to first order in $\hbar$, which makes
$C_1$ be a closed Hochschild cochain, namely
$$
\partial C_1(f,g,h) := fC_1(g,h) - C_1(fg,h) + C_1(f,gh) - C_1(f,g)h = 0.
$$
If $C_1$ is exact; that is, if there is a 1-cochain such that
$$
C_1(f,g) = \partial E(f,g),
$$
then to first order in $\hbar$ Eq.(2.4) can be written
$$
(f - \hbar E(f))*(g - \hbar E(g)) = fg - \hbar E(fg);  
$$
essential first order deformations
are classified by Hoch$^2$. [G1]

Suppose that a formal  $*$-product is associative to order $p \geq
1$; this statement involves $C_1,\cdots ,C_p$ only, and we suppose these
cochains fixed. Then the condition that must be satisfied by
$C_{p+1}$, in order that the $*$-product be associative to order $p+1$, is
$$
\sum_{m,n = 1\atop m+n = p+1}^p  \biggl(
C_m(C_n(f,g),h) -C_m(f,C_n(g,h)\biggr) = \partial C_{p+1}(f,g,h).\eqno(2.5)
$$
The left hand side is closed [G1],[BFFLS]; this minor miracle is responsible
for the fact that commutative algebras are not isolated in the family of
associative algebras. The right hand side of Eq. (2.5) is a Hochschild
coboundary. An obstruction to the existence of a two-cochain
$C_{p+1}$ that  would solve Eq. (2.5) is thus an element of  Hoch$^3$.  
This statement will be sharpened below.  
 
\ve 

\ce{\bf  2.2. 
The BGS decomposition of Hochschild (co-)homology.}

The $p$-chains of the Hochschild homology complex of a commutative algebra $A$
are the $p$-tuples $a = \sum a_1\otimes \cdots \otimes a_p \in A^{\otimes
p}$, and the
differential is defined by
$$
da = a_1a_2\otimes a_3 \otimes \cdots \otimes a_p - a_1\otimes
a_2a_3\otimes a_4\cdots \otimes a_p + \cdots +(-)^pa_1\otimes \cdots
a_{p-2}\otimes a_{p-1}a_p.
$$
The $p$-cochains are maps $A^{\otimes p}   \rightarrow  A$,  and the
differential is
$$
\partial C(a_1,\cdots ,a_{p+1}) = a_1C(a_2,\cdots ,a_{p-1}) - C(da) +
(-)^{p+1}C(a_1,\cdots ,a_p)a_{p+1}.
$$

The
Hochschild cochain complex splits into a finite or infinite sum of direct
summands. (If the algebra is generated by $N$ generators then there are only
$N$ nonzero summands.)     After the pioneering work of Harrison [H] and Barr
[B],[G2], the complete decomposition of the Hochschild cohomology of a
commutative algebra was found by Gerstenhaber and Schack [GS]. 
 The decomposition is based on the action of $S_n$
on
$n$-cochains, and on the existence of  $n$ idempotents
$e_n(k),~k = 1,\cdots n,$~in $\Crm S_n$,   $\sum_ke_n(k) = 1$,  with
the property that
$$
\partial\circ e_n(k) = e_{n+1}(k)\circ \partial.
$$
There is thus a decomposition ${\cal C}^n = \sum_{k = 1}^n {\cal C}^{n,k}$ of
the space  of $n$-cochains,  and Hoch$^n = \sum_{k = 1}^n {\cal H}^{n,k}$ 
with   ${\cal H}^{n,1} =:~$Harr$^n$.

A  generating function was found by Garsia [G],
$$
\sum_{k = 1}^nx^ke_n(k) = {1\over n!}\sum_{\sigma\in S_n}
(x-d_\sigma)  (x-d_\sigma +1)  \cdots   (x-d_\sigma +n-1){\rm sgn}(
\sigma)\sigma,
$$
where $d_\sigma $ is the number of descents,~~$\sigma(i) > \sigma(i+1)$
   ,~~in $\sigma (1 \cdots n)$. (Example: $\sigma(1234)
= 3142$ has one descent, from 2 to 3.) The simplest idempotents are
$$\eqalign{ 
e_2(1) 12  &=  {1\over 2}(12 + 21),~~e_2(2)12 = {1\over 2}(12 - 21),\cr  
e_3(1) 123&=  {1\over 6}\bigl(2(123 - 321)  + 132 - 231 + 213 - 312
\bigr), \cr e_3(2) 123&=  {1\over 2}(123 + 321)\cr
e_n(n)& =  {1\over n!}\sum_{\sigma\in S_n}{\rm sgn}( \sigma)\sigma.
\cr}
$$
The space of Hochschild $n$-chains decomposes in the same way, ${\cal C}_n =
\sum_{k = 1}^n {\cal C}_{n,k}$ with 
$$
d\circ e_n(k) = e_{n-1}(k)\circ d
$$ 
so that Hoch$_n = \sum_{k = 1}^n{ \cal H}_{n,k}$  with  ${\cal H}_{n,1} =
$Harr$_n$.

\b\b 

\ce{\bf 2.3. Star products and the BGS decomposition.}  

 Every 2-cochain $C$ has a decomposition
$$
C = C^+ + C^-, ~~ C^+ \in\, {\cal C}^{2,1},~ C^- \in {\cal C}^{2,2}.
$$
Associativity of the star product to  order $\hbar$ is the requirement that
both 2-cochains be closed, $\partial C^+ = \partial C^- = 0$.  A first order
deformation of the algebra is inessential if both forms are exact.  In the
case of a smooth manifold Harr$^2$ is empty and choosing $C_1^+ = 0$ entails no
essential loss. Returning now to   Eq.(2.5) we have  mappings
$$
\hskip-2cm\hskip-1cm {\cal Z}^{2,1}\times {\cal Z}^{2,1}
$$
\vskip-1.4cm
$$
\hskip-2cm \hskip.6in \searrow
$$
\vskip-1.cm
$$  
\hskip-2cm\hskip1.9in{\cal Z}^{3,1} ~ ~ \leftarrow ~  {\cal Z}^{2,1} 
$$
\vskip-.8cm
$$ 
\hskip-2cm \hskip-.4cm{\cal Z}^{2,2}\times {\cal Z}^{2,2} ~~   {\nearrow\atop
\searrow}\eqno(2.6)
$$
\vskip-.7cm
$$   
\hskip-2cm\hskip1.27in{\cal Z}^{3,3} 
$$
\vskip-.8cm
$$
\hskip-2cm\hskip2.1 cm {\cal Z}^{2,1}\times {\cal Z}^{2,2} ~~ \rightarrow ~ 
 {\cal Z}^{3,2} ~~ \leftarrow ~   {\cal Z}^{2,2}  
$$                
 The first column of arrows represents the construction  on the left side of
Eq. (2.5).  The second column of arrows is the mapping by the differential.
Because
${\cal B}^{3,3}$ is empty; the obstruction in ${\cal Z}^{3,3} = {\cal H}^{3,3}$
 demands that the antisymmetric part of the left side of (2.5) 
vanish, this  is
the Jacobi identity, satisfied if $C_1^-$ is a Poisson bracket.   The
first line shows that the obstruction to abelian deformations is the Harrison
component Harr$^3 = {\cal H}^{3,1}\subset ~$Hoch$^3$. In the case of smooth
manifolds Hoch$^n = {\cal H}^{n,n} = 0$ and abelian deformations are
inessential. The familiar  deformations with $C_1 = C_1^-$
encounter no additional obstructions to order $\hbar^2$.  We now
turn to a preliminary investigation of  algebraic varieties with singularities.
We shall find varieties for which ${\cal Z}^{3,1}$,~${\cal Z}^{3,2}$ and
${\cal Z}^{3,1}$ are all non empty.
\ve

\no{\steptwo  3. Some varieties with singularities.}~
 
\ce{\bf  3.1. Conic varieties defined by quadratic relations.}~ 
These are algebraic varieties of the type~ $\Cit^N/R$, where $R =
\{g_\alpha\}_{\alpha = 1,2,...}$ is a set of homogeneous, quadratic forms,
$$
g_\alpha = \sum_{i,j = 1}^N g_\alpha^{ij}x_ix_j,~~\alpha = 1,2,...~.\eqno(3.1)
$$
Let $A$ be
the graded coordinate  algebra $A  = \Cit[x_1,...,x_N]/(R)$, and  
 $A_+$ the subalgebra of positive degrees; the restriction to positive
degrees is essential.   Cochains on $A_+$ extend naturally to $A$, but the
homology of
$A_+$ is richer than that of $A$. (The generators of $ A_+$ are not exact.)

\b
\no  {\bf 3.1.1. Definition.}~{\it   A chain $a = a_1\otimes
a_2\otimes ...\otimes a_p$ will be said to be `linear' if each $a_k, k =
1,...,p$ is of degree 1. The `restricted complex' is  the restriction of the
Hochschild complex of
$A_+$ to closed, linear chains.}
\b 
In the case of a smooth manifold, or more generally in the case of a regular
commutative algebra, the chains of the restricted complex are the
antisymmetric  ones;  the next result  reduces in that case to a
famous theorem of Hochschild, Kostant and Rosenfeld [HKR]. For  other
generalizations see [FG],[FK] and Sect. 3.3.

Let $A_+$ be as above, the subalgebra of $A = \Cit[x_1,...,x_N]/(R)$ obtained
by restriction to positive degrees.
\b
\no {\bf 3.1.2. Theorem.}{\it   ~ The Hochschild complex of $A_+$
is quasi-isomorphic to the restricted complex; that is, their
(co-)homologies are isomorphic.}  
\b
\no{\bf  Proof.}~
The restriction of a
closed/exact form is closed/exact. Conversely, every closed/exact restricted
cochain is the restriction of a closed/exact Hochschild cochain. It
is enough to consider homogeneous chains; that is,  $a = a_1\otimes ...\otimes
a_n$ such that each factor $a_k,\, k = 1,...,n$  is of well defined degree.
(The only grading that we use is the total polynomial degree.) To show that
every restricted  (= closed, linear)
$n$-cochain extends to a closed, Hochschild $n$-cochain we consider the formula
$$
\partial C(a_1,...,a_{n+1}) = a_1C(a^1) - C(da)
+(-)^{n+1}C(a^{n+1})a_{n+1},\eqno(3.2)
$$
where $a^1 = a_2\otimes ... \otimes a_{n+1}$ and $a^{n+1} = a_1\otimes
...\otimes a_n$. Evidently the degree of $a$ is higher
than the degrees of
$a^1$ and $a^{n+1}$.  
This formula can therefore be used to try to extend closedness,
recursively  to higher degrees. The obstruction is
$da = 0$, but it is easy to verify that the remaining terms in (3.2) vanish
when
$a$ is exact and $C$ is closed on lower degrees. The  obstruction comes
from homology; if the restricted $n$-cochain $C$ has the
property that
$\partial C(h) = 0$ for a representative
$h$ of every homology class of $n+1$-chains, then it extends to a closed,
Hochschild  $n$-cochain.
 \b

\ce{\bf  3.2. The simple cone; one quadratic relation.}~ Retain all
the definitions but suppose that $R = g$ is just one quadratic form. In this
case
 the following holds.
\b
\no{\bf  3.2.1.  Proposition.}  ~  {\it  Every closed chain is
homologous to a linear chain and no linear chain is exact. The space ${\cal
Z}_{2k+l}$   of closed, linear $(2k+l)$-chains is spanned by the following
$n = 2k+l$-chains, $n$ = 1,2,...~,} with $Z_{2k+l,k+l} \in {\cal
Z}_{2k+1,k+l}$,
$$\eqalign{&
(Z_{2k+l,k+l})_{m_1...m_l} =
g^{i_1j_1}...g^{i_kj_k}\sum_\sigma(-)^\sigma x_{i_1}\otimes ...\otimes
x_{j_k}\otimes x_{m_1}\otimes ...\otimes x_{m_l},\cr &
\hskip1in m_1,...,m_l= 1,...,N,\cr}
$$      
{\it where the sum is over all permutations of $i_1...j_km_1...m_l$ that
preserve the internal order of each pair 
$(i_1,j_1), (i_2,j_2),...,(i_k,j_k)$.}
\b
 Examples,
$$
(Z_{1,1})_i= x_i,~~ i = 1,...,N,
$$
$$
(Z_{2,2})_{ij} =  x_i\wedge x_j,~~ Z_{2,1}=~ g^{ij}x_i\otimes x_j,
$$
$$
(Z_{3,3})_{ijk} = x_i\wedge x_j\wedge x_k,~~  
(Z_{3,2})_{k} =g^{ij}(x_i\otimes x_j \otimes x_k - x_i\otimes x_k \otimes x_j 
+ x_k\otimes x_i \otimes x_j ), 
$$
$$
(Z_{4,4})_{ijkl} =  x_i\wedge x_j\wedge x_k\wedge x_l,
$$
$$\eqalign{ 
(Z_{4,3})_{kl} =~& 
g^{ij}( x_i\otimes  x_j\otimes x_k\otimes x_l +
x_k\otimes x_i\otimes x_j\otimes x_l + x_i\otimes x_l\otimes x_j\otimes
x_k\cr & 
 + x_k\otimes x_l\otimes x_i\otimes x_j + x_i\otimes x_k\otimes
x_l\otimes x_j + x_l\otimes x_i\otimes x_k\otimes x_j- k,l),
\cr}
$$
$$ 
(Z_{4,2})= g^{ij}g^{kl}(x_i\otimes x_k\otimes x_l\otimes x_j -
x_i\otimes x_k\otimes x_j\otimes x_l + x_i\otimes x_j\otimes x_k\otimes x_l).
$$

 From this point onwards, the notation  ${\cal Z}_{k,l}, {\cal B}_{k,l}, {\cal
Z}^{k,l}, {\cal B}^{k,l}$ stands for spaces defined with reference to the
restricted complex. The dimension of
${\cal Z}_{2k+l,k+l} \subset {\cal H}_{2k+l,k+l}$ is
$\Big({N\atop l}\Big)$ or 0. We note that ${\cal Z}_{3,1}$ is empty.
To lowest order, a star product is determined by the 2-cochain $C_1$, and up to
equivalence  by the restricted 2-cochains; that
is, by their values $C_1^-(x_i\wedge x_j)$ and
$g^{ij}C_1^+(x_i, x_j)$ on the homological basis.  The differential
$\partial C_1^+$ is in
${\cal B}^{3,1}$, and since ${\cal Z}_{3,1}$ is empty every restricted,
symmetric 2-cochain is closed.   The
differential
$\partial C_1^-$ is in
${\cal B}^{3,2}$ ; $C_1^-$  is closed iff $\partial C_1^-$  vanishes on
$Z_{3,2}$, 
$$\eqalign{&
 \partial C_1^-(Z_{3,2}) =  4g^{ij}  x_iC_1^-(x_j, x_k)  = 0.\cr}  
$$
This will be interpreted as the statement
that the `Hamiltonian' vector fields $C_1^-(x_k,\cdot)$ be tangential to the
constraint surface.  

The addition of an exact form
$\partial  E$ to $C_1$ (we like to think of it as a `gauge
transformation') does not affect $C_1^-$ but it adds
$2g^{ij} x_iE(x_j)$ to $C^+_1(Z_{2,1}) = g^{ij}C_1^+(x_i,x_j)$.  
\b
\no{\bf  3.2.2. Example.}~Suppose that the
2-form $g$ is nondegenerate; then $C_1^+$  is fixed up to equivalence
by its value
$c = g^{ij}C_1^+(x_i,x_j)|_{x = 0}$ at
$x_1 = x_2 = ...  = 0$. That is;
${\cal H}^{2,2}$ is the space of tangential vector fields on the cone and
${\cal H}^{2,1} = ~\Cit$. (See 3.2.8 for equivariant cohomology.) 
 \b
We examine the obstructions to extending the star product to all orders in
$\hbar$, referring to Eq.s~(2.5) and (2.6).
\b
\no{\bf  3.2.3.}~
 The emptiness of ${\cal B}^{3,3}$ is an
obstruction that must be circumvented by imposing the Jacobi identity on
$C_1^-$.  Because $C_1^-$ is closed this entails that it extends to a
unique Poisson bracket on $A$. Recall that, if $(a,b) \mapsto \{a,b\}$ is a
Poisson bracket on $A$ then $\forall f \in A$ the mapping $f^\#: A\rightarrow
A$   defined by $a \mapsto \{f,a\}$ is a derivation.
The projection of (2.5) on ${\cal Z}_{2,3}$ is now solved by taking
 $
C_n^-(x_i,x_j) = 0,~~i,j = 1,...,N,~ n>1. 
 $
This choice  is implicit in the context of differentiable deformations,
and it is one of the axioms of invariant star products. We do
not investigate alternatives.

\b
\no Summary: The projection of Eq. (2.6) on ${\cal Z}_{3,3}$ leads to

(a) $C_1^-(x_i,x_j) = \{x_i,x_j\}$ ~ extends  to  a  Poisson bracket on $A$,

(b) $\{x_i,g(x,x)\} = 0$\hskip1.1cm is the condition that $C_1$ be closed,
 
 (c)  $C_n^-(x_i,x_j) = 0,~~ n > 1,~~ i,j = 1,...,N$.

\vskip-1.2cm $$\eqno(3.3)$$

 \b 
\no{\bf  3.2.4.} ~ The projection of Eq. (2.5) on ${\cal Z}_{3,2}$
takes the form
$$\eqalign{ 
\sum_{m+n = p+1}\Big[-4g^{ij}& C^+_m(C^-_n(x_i,x_k),x_j)   +
2g^{ij}C_m^-(C_n^+(x_i,x_j),x_k)\Big]\cr &    = \partial
C^-_{p+1}(Z_{3,2})_k   = 4g^{ij}x_iC_{p+1}^-(x_j,x_k) ,~~k  =
1,...,N.\cr}\eqno(3.4)
$$
If $g$  is nondegenerate then we can restrict the value of $C_1^+(x_i,x_j)$
to \Crm; see Example 3.2.2. The obstruction is then the value of the left side
at $x = 0$. In view of (3.3), (3.4) simplifies to
$$
 2g^{ij}  C^+_p( \{x_i,x_k\},x_j) - g^{ij}\{C_p^+(x_i,x_j),x_k\} 
 = 0.\eqno(3.5) 
$$
 \b
\no{\bf 3.2.5. Example.}~ If $\{x_i,x_k\} = {\epsilon_{ik}}^m$, the
coefficients $({\epsilon_{ik}}^m)$ the structure constants of a simple Lie
algebra \g, then (3.5) is satisfied when $C^+_p$ is the Killing form of \g.
 \b
\no{\bf 3.2.6. Example.}~ In the case of sl(n) there is an
equivariant 3-tensor $f: A\otimes A \rightarrow A$ and $C^+_p = f$ also solves
(3.5). 
\b
\no{\bf 3.2.7. Example.} ~ 
 ~  Choose
coordinates such that $g(x,x) = x_N^2 - \rho(x)$, where $\rho(x)$ is a
polynomial in $x_1,...,x_N$, at most linear in $x_N$.  
   A
regular function on the cone, being the restriction of a polynomial in
$x_1,...,x_N$, has a unique decomposition of the form
$f = f_1 + x_Nf_2$, where $f_1,f_2$ are polynomials in $x_1,...,x_{N-1}$.
Define ([F2]) 
$$
f*g = fg + \hbar f_2g_2 = fg|_{x_N^2 =\rho(x) + \hbar}.\eqno(3.6)
$$ 
This star product is associative to all orders; it can be interpreted as  
the ordinary\break  product of functions on the hyperboloid $g = \hbar$. This 
deformation, in which $C_1$ is symmetric, can be followed by another
deformation in which
$C_1$ is antisymmetric, leading to a Poisson bracket such that the vector
fields $\{x_i,\cdot\}$  are tangent to the hyperboloid
$g(x,x)=
\hbar$;  no Harrison cohomology intervenes in either stage.   When both
deformations are combined we note that, with
$C_1^+$ as we have defined it in Eq.(3.6), $C_1^+(f,g) = f_2g_2$, the only
contribution to the left side of Eq.(3.5) at $x = 0$ comes from the linear term
in
$\{x_i,x_k\}$.
Closure of $C_m^-$ implies that the vector field $\{\cdot,x_k\}$
is tangent to   $g(x,x) = \hbar$.

 \b
\no{\bf  3.2.8. Equivariant cohomology.}~ In the context of Lie
algebras and invariant star products all maps will be equivariant. This
affects the question of exactness, as in 3.2.2, and will be taken into account
as the occasion arises. Theorem 3.1.2 is not affected.   
\b\b

\ce{\bf  3.3. Case of one polynomial relation.}

Let $A = \Cit[x_1,...,x_N]/(g)$ where
$g$ is a polynomial without constant or linear terms.
\b

\no{\bf  3.3.1. Grading.} ~We can no
longer restrict our attention to homogeneous chains.   By an
appropriate linear transformation of variables we can bring the polynomial $g$
to the form
$x_1^t +  h(x_1,...,x_N)$, where
$h$ is a polynomial of degree less than $t$ in $x_1$. A `normalized' polynomial
is one that is of degree less than $t$ in $x_1$; to each regular function on
$\Rit/g$ there is just one normalized polynomial and the degree of a regular
function is defined to be the degree of this normalized polynomial.

Let $g$ be a normalized polynomial in $x_1,...,x_n$, without a constant term
and without linear terms. Choose a presentation
$$
g = \sum_{a,b = 1}^Kg^{ab}y_ay_b,~~g^{ab} = g^{ba}\in ~\Cit,~ a,b = 1,...,K. 
$$
where each of the polynomial factors $y_a, ~a = 1,...,k$,
has no constant term, and let $A$ be the filtered algebra $A =
\Crm[x_1,...,x_N]/(g)$ with $A_+$ the sum of the positive degrees of $A$.
To each element of $A$ is associated a unique, normalized polynomial.

 \b
 Let
${\cal Z}_{2k+l}$  denote the space of closed  $p$-chains,  spanned by the
following  
$p$-chains, $p = 2k + l  = 1,2,...$ and $~ m_1,...,m_l= 1,...,N$,    
$$\eqalign{& 
\hskip-1cm  (Z_{2k+l,k+l})_{m_1...m_l}  =
g^{a_1b_1}...g^{a_kb_k}\sum_\sigma(-)^\sigma \sigma(y_{a_1} 
 \otimes ...\otimes  y_{b_k}\otimes x_{m_1}\otimes 
...\otimes x_{m_l}), 
  \cr} 
$$        
 where the sum is over all permutations  
that preserve the internal order of each pair    
$(a_1,b_1),...,(a_k,b_k)$.   
\b
\no{\bf   3.3.2.  Theorem.}    ~ {\it The restriction of the
Hochschild complex to such chains is a quasi homomorphism.}  
\b

\no The proof of the theorem is as the proof of 3.1.2, but  needs  the
following lemma.
 Let  $A_+$ be the filtered algebra 
without unit as above. The degree of an $A_+$-chain $ a_1\otimes ...\otimes
a_p$ is the sum of the degrees of its factors. 
\b
\no{\bf  3.3.3. Lemma.}~ {\it If an
$A_+$ $p$-chain    is exact, then there is a $p+1$-chain $b$,} with the same
degree as $a$, {\it such that $a = db$.} [FK]
\b\b

\no{\steptwo  4. Invariant star products on coadjoint orbits.} 

\ce{\bf  4.1. Background.} ~The origin of this problem is as follows
[AM]. Let
$W$ be a symplectic space with Poisson tensor $\Lambda$
and let there be given an action by so(3), generated by  
hamiltonian vector fields $\Lambda(dL_i)$, where $L_1,L_2,L_3 \in C^\infty
W$  satisfy the following Poisson bracket relations,
$$
\{L_i,L_j\} = \epsilon_{ijk} L_k.
$$ 
The Casimir element
$$ 
Q = \sum_i(L_i)^2
$$ is an invariant; that is,  $\{L_i,Q\} = 0$. The hamiltonian vector fields
leave invariant each surface $Q =  ~$constant, each such surface is a
symplectic leaf with an induced Poisson structure. The problem: to invent an
equivariant ordering   (an invariant star product) such that the star
polynomial ${\cal W}(Q) := \sum L_i*L_i$ is invariant, $L_i*{\cal W}(Q) =
{\cal W}(Q)*L_i$, and fixed,
$\sum
L_i*L_i*f = qf$ for some $q\in \Cit, \forall f$. 
\b

\ce{\bf 4.2. Invariant star products.}

  Let $G$ be a Lie group,  \g ~ 
the  Lie algebra of $G$,  \g$'$  the real vector space dual, and
$W$ an orbit of the coadjoint action of the connected component of $G$ in 
 \g$\,'$. This
  defines a   homomorphism from the symmetric algebra
$S(\hskip-0.2mm\g)$  into $C^\infty W$. There is a natural Poisson structure
on $W$, such that $\hbar \{a,b\} = [a,b]$ for $a,b\in  \,$\g. (We identify
\g~ with \g$'$.)

\b
\no{\bf  4.2.1.  Definition.} ~  {\it A star product on a
  coadjoint orbit  
$W$ is   \g -invariant if, for all} $k\in \,\Cit;\,$ $ \it a,b \in\, $\g$\,;
\,f,g
\in C^\infty W$, 
$$
\eqalign{&
k*a = a*k = ka,\cr
&
a*b - b*a =  \hbar \{a,b\},\cr
&
\{a, f*g\} = \{a,f\}*g + f*\{a,g\}.
\cr}\eqno(4.1)
$$
\b

\no{ \bf 4.2.2. Remark.}   ~Coadjoint orbits provide a plethora of
symplectic spaces, but to invoke the assistance of a Lie group \underbar {for
that purpose alone} is somewhat odd. It seems more natural, in this
context, to investigate star products that incorporate additional elements of 
group theoretical structure. 
\b 
Given an associative star product on $W$,
a linear map ${\cal W}$  from the symmetric algebra into $C^\infty W$ is
defined as follows, 
$$
{\cal W}: a^n \mapsto {\cal W}(a^n) = a*a*... \in C^\infty W,~~ a\in  \g,
$$
Conversely, any  invertible linear map ${\cal W}$ that associates a
$C^\infty$  function to each formal star monomial defines an associative 
star product on
$C^\infty W$. For any polynomial $P(a)$ we write $P(a,*)$ for ${\cal
W}(P(a))$.

 The original Moyal product is  the unique associative, invariant product 
for $W = \Rit^{2N}$  with the standard Poisson bracket  such that
$ {\cal W}(a^n) = a^n$ for every $a$ that is linear in the natural
coordinates.   It is invariant under the Lie algebra of affine symplectic
transformations. The domain includes the space of regular functions (the space
of polynomials in the generators).

 A recipe for the construction of all invariant star products for
any compact, semisimple Lie algebra, on any regular coadjoint orbit, was
formulated  almost 30 years ago. (For non-regular orbits of a compact group
see [Ll].)
\b

\b
\no{\bf 4.2.3. Definition.}{\it   ~~A star product on $W$ is
nondegenerate if the space of star polynomials (actally, the image by ${\cal
W}$) is dense in the space of
$C^\infty$ functions on $W$.}

\b
\no{\bf  4.2.4. Theorem.} ~{\it An  associative,
nondegenerate, invariant star product on a coadjoint  orbit $W$ of a compact,
semisimple Lie algebra
 \g~  is given by  an infinitesimal character $Z(\g) \mapsto
\Cit$ and  the
 formulas
$$\eqalign{& 
P_n(a,*) :=  C_nP_n(a), ~~ C_n\in \,\Cit-\{0\},~~  C_0 = C_1 =
1, \cr & a*b - b*a = \hbar[a,b],~~  a,b\in\g,~~n =
0,1, 
\cr}\eqno(4.2)
$$
where $P_0 = 1,~ P_1(a) = a,~\{P_2,P_3,...\}$ is a complete set of 
irreducible, harmonic elements of  $S(\g\, ')$ and $P_n(a,*)$ are the
corresponding star polynomials. }    [BFFLS]
\b
In the case when  \g = so(3), $\, P_n(a)$ is a solid Legendre polynomial and
the polynomials $P_n(a,*)$ can be obtained from the recursion relation
(found and solved in [BF]) 
$$
(n+1)P_{n+1}(a,*) = (2n+1)\,a*P_n(a,*) - n(q +{1-n^2\over 4}
 \hbar^2)|a|^2P_{n-1}(a,*), 
$$
with $P_0(a,*) = 1$. The parameter $q$ is the value of the Casimir operator
$\sum L_i*L_i$.

 The statement of the theorem remains valid for
noncompact, semisimple Lie algebras and regular orbits. However, the
restriction to polynomials would not be  appropriate.

\b
\ce{\bf  4.3. Finite dimensional representations.}

An invariant star product gives an action of $\g\,$ on the star algebra, and on
$C^\infty(W)$, by the homomorphisms $\pi_l: a\mapsto a*$ and similarly by
$\pi_r: a\mapsto  *(-a)$,
$$
\pi_l(a)f = a*f,~~ \pi_l(a)f = -f*a,~~ f \in C^\infty W.
$$
In the case of compact Lie algebras we expect to find finite dimensional
representations. How this actually comes about can  be seen explicitly 
in the case $\g = so(3)$.
\b  
\no{\bf  4.3.1.   Example. Proposition.}  
{\it Let $\g = so(3)$; an invariant, associative star product is of
one of two types, both defined as in Eq.(4.2). 

(1) All $C_n \neq 0$; the
action 
$\pi_l $ or $\pi_r$ generated by
$a*~$or$  ~  *(-a), a \in \g$ is not semisimple. In this case the choice
 $C_n = 1 + o(\hbar)$ for   $n>1$ provides an equivariant deformation
for every value $q$ of the Casimir $\sum L_i*L_i$. 

(2) If the infinitesimal character takes the Casimir to the value $q = l(l+1)$,
$2l \in \{0,1,2,...\}$  
then the algebra of star polynomials contains an ideal generated by
$P_{2l+1}(a,*)$. The quotient is a finite dimensional $*$-algebra and is
spanned by the projection of $\{P_n(a,*)\}_{n = 0,1,...,2l}$. The action of
$a*$ and $ *(-a)$ is equivalent to the direct product of two copies of the
irreducible representation of $su(2)$, each with dimension $2l+1$. The
polynomial
$P_{2l+1}(a,*)$ reduces in this case to 
$$
 P_{2l+1}(a,*)\propto \prod_{m = -l}^l (a*-m|a|),~~ a = \sum a_iL_i,~ |a|:=
 \sqrt{\sum a_i^2},
$$
and every $P_n(a,*)$ with $n > 2l+1$  contains $P_{2l+1}(a,*)$ as a factor. }
\b
The spectrum, the space of maximal ideals, has a
finite number of disconnected components.  
\ve
\no{ \steptwo 5. Introduction to singular orbits.} 
 
\ce{\bf  5.1. The coadjoint orbits of so(2,1).}

Let $L_1,L_2,L_3$ be the standard basis for the real Lie algebra $so(2,1)$,
with relations $[L_i,L_j] = \epsilon_{ijk}L_k$ and Killing form
$g^{ij}L_iL_j = -(L_i)^2 -(L_2)^2 + (L_3)^2$.    
The moment map interprets $L_1,L_2,L_3$ as coordinates for the coadjoint
space, and in this role we denote them by the symbols $x = (x_1,x_2,x_3)$.
The regular orbits are the loci of $g(x) = g^{ij}x_ix_j =~ $const$~ \neq 0$ and
the nilpotent orbit  is the algebraic variety
$$
M = \Rit^3/g(x).
$$
A regular function on $M$ is the restriction of a polynomial in
$x_1,x_2,x_3$. This variety is a simple cone in the sense of Section  3.2. The
first order star product $x_i*x_j = x_ix_j + \hbar C_1(x_i,x_j)$ is equivariant
if
$C_1^-(x_i,x_j) = (1/2)
\epsilon_{ijk}x_k$ and
$C_1^+(x_i,x_j) \propto  g_{ij}$. Since $g$ is nondegenerate there is a two
parameter family of essential, first order equivariant deformations, indexed by
$\hbar$ and $C_1^+|_{x = 0}$. Every first order, equivariant deformation can be
extended to an invariant star product (to all orders in $\hbar$); for
example,   by the method outlined above. Questions of domains have not yet
been adequately discussed, to our knowledge.

\b\b
 \b

\ce{\bf  5.2.  Minimal orbits of sl($n$),$~ n > 2$.}~  If the
matrix
$U\in
\,$sl($n$)  lies on the minimal orbit, then $U^2 = 0$, but this relation is
not enough to define a minimal orbit.   The minimal orbit is defined by the
relations 
$$
U_a^bU_c^d - U_a^dU_c^b = 0,~ a,b,c,d = 1,...,n.\eqno(5.1)
$$
These are solved by the {\it factorization} 
$$
U_a^b = p_aq^b,~~ {\rm with}~ ~q\cdot p := q^ap_a = 0,\eqno(5.2)
$$
which defines an imbedding of $\g$ into  the space of second
order polynomials on $P^{2n-2} =\Rrm^{2n}/(p\cdot q, \approx)$, where $\approx$
is the equivalence relation $\forall\lambda > 0: (\lambda p,\lambda^{-1}q)
\approx (p,q)$. A star product can be defined on this space by introducing the
Poisson bracket defined by $\{q^b,p_a\} = \delta^b_a,~\{q^b,q^d\} =
\{p_a,p_c\} = 0$ and quantizing this in the manner of Weyl, see Sections 5.4
and 6.3.  
\b

\ce{\bf  5.3 Minimal orbits of sp($2n$), so($n$) and the others.} 

Sp(2$n$) is the algebra of traceless matrices that leaves  invariant an anti
symmetric, nondegenerate 2-form $\eta\,$:  $ U_a^b\,\eta_{bc} =: L_{ac} =
L_{ca}$. The imbedding $L_{ac} = \xi_a\xi_c$, with $ \xi_1,...,\xi_{2n}  =
q^1,...,q^n,p_1,...,p_n $ incorporates all the
relations that define the minimal orbit, namely
$$
L_{ab}L_{cd} = L_{ad}L_{cb}.
$$

The Lie algebra so($n$) is as  sp(2$n$)  except that the form $\eta$ is
symmetric. Analogy with sp(2$n$) suggests using Grassmann
variables, replacing the commutative affine algebra by a supercommutative
super Lie algebra, as is done in field theories with fermions. An alternative
is the imbedding
$L_{ac} = q_ap_c - a,c$ (with $q_a = q^b\eta_{ba}$). The relations that define
the minimal orbit are  
$$
\sum_{cyc(abc)}L_{ab}L_{cd} = 0,~~ \eta^{bc}L_{ab}L_{cd} = 0.\eqno(5.3) 
$$
The first is implied by the embedding and the second can be incorporated by
restriction to $\eta^{ab}p_ap_b = \eta^{ab}q_aq_b = q^ap_a = 0$ and
projecting on a quotient. 

The Lie algebra  G$_2$ is a subalgebra of so(7). The minimal orbit can be
parameterized as that of so(7) with the additional condition $p\times q = 0$.
 
\b
 
\ce {\bf  5.4. Associated representations and Joseph ideals.} 

\no{\bf  5.4.1. Background.} ~The ideals in the enveloping algebra
of a compact, simple Lie algebra are fixed by a central character;
the noncompact
case is more interesting. An example is well known to physicists.
Consider the Lie algebra
$so(4,2)$, with the usual basis $\{L_{ab}\}_{a,b = 1,...,6}$ and relations
$$
[L_{ab},L_{cd}] = \eta_{bc}L_{ad} - a,b - c,d~,\eqno(5.4)
$$
where $\eta$ is the pseudo Euklidean metric. The tensor $g^{bc}(L_{ab}L_{cd} +
a,d)
$ in the enveloping algebra reduces, in a certain irreducible and unitarizable 
representation, to fixed numerical values, so that  relations of the type
 $$
\eta^{bc}(L_{ab}L_{cd} + a,d) = -\hbar^2 \eta_{ad} \eqno(5.5)
$$
hold  in the representation. This particular representation is well known from
its appearance in the analysis of the Schroedinger theory of the hydrogen atom
and in conformal field theory.
The algebra also enters the description of Keplerian orbits on a
6-dimensional phase space.  

This orbit is of interest, {\it inter alia}, because Kostant's method of
geometric quantization encounters a difficulty, the non existence of an
invariant quantization [K1]. (It was shown by Joseph that this is true
of all minimal orbits except the case of sl(n) [J].)  Although the
corresponding quantum theory is known, an invariant Wigner-Weyl correspondence
is not. The existence of an invariant star product associated with such a
correspondence is strongly expected to exist, but it has not been constructed.
Nevertheless,   the relation (5.5) suggests that there is an invariant star
product such that
$$
\eta^{bc}(L_{ab}*L_{cd} + a,d) = -\hbar^2 \eta_{ad}. \eqno(5.6)
$$
A question that motivated this work is whether such a deformation  exists, and
if it is a deformation in the direction of the Poisson bracket. The
undoubted presence of interesting homology on this highly singular orbit was
expected to play a role in invariant quantization. A preliminary exploration 
of the associated representations will show us what to expect.

\b

\def\d{\delta}

\no{\bf  5.4.2. Associated representations of $sl(n)$.}
 
Let $V_N$ denote a space of functions on $\Crm^n$, spanned by a set of
functions $x_1^{r_1},...,x_1^{r_n}$ with $N = \sum_i r_i$ fixed, and let
$\{\tilde U_a^b\}_{a,b = 1,...,n}\}$ be the family of operators given by
$\tilde U_a^b =
\hbar(x_a\partial/\partial x_b - (N/n)\delta_a^b)$ in $V_N$. Then
$[\tilde U_a^b,\tilde U_c^d] =
\hbar(\delta_c^b\tilde U_a^d - \delta_a^d\tilde U_c^b)$, which are the
relations of
$sl(n)$ with the association that identifies $\tilde U_a^b/\hbar$ with the unit
matrix
$E_{ab}$ when $a
\neq b$ and with $E_{aa} - 1/N$ when $a = b$.

To be more precise, consider the real form $su(n-1,1)$, with the compact
subalgebra $su(n-1)$ generated by $\{U_a^b\}_{a,b = 1,..., n-1}$. Taking
$ r_1,...,r_{n-1}$ to run over the natural numbers one obtains, for a range of
values of $N$, a unitarizable, highest weight representation, with the highest
weight reducing to zero on $su(n-1)$. This representation is finite dimensional
if
$N
\in
\{0,1,... \}$. The relations that define the minimal orbit are $U_a^bU_c^d -
b,d = 0$, while
$$
\tilde U_a^b\tilde U_c^d -b,d =   \hbar\d_c^b\tilde U_a^d - \hbar{N\over
n}(\d_a^b\tilde U_c^d + \d_c^d\tilde U_a^b) - \hbar^2{N\over
n}({N\over n} + 1)\d_a^b\d_c^d  - b,d.
 \eqno(5.7)
$$
An invariant deformation in the direction of the Poisson bracket would have
$$
U_a^b*U_c^d = U_a^bU_c^d +(\hbar/2)\{U_a^b,U_c^d\} + \hbar C^+(U_a^b,U_c^d).
 \eqno(5.8)
$$
In this setting, because of the very strong relations that characterize the
orbit,  the most general equivariant, symmetric 2-form $C^+$ takes the form
$$
\hbar C^+(U_a^b,U_c^d) = (k/2)\big(\delta_c^bU_a^d + \delta_a^dU_c^b-{1\over
n}(\delta_a^bU_c^d + \delta_c^dU_a^b)\big) + k'(\delta_a^d\delta_c^b - {1\over
n}\delta_a^b\delta_c^d).\eqno(5.9)
$$
This yields relations just like (5.7) if the parameters $k,k'$ 
are appropriately related to the degree $N$ of the homogeneous
functions in the vector space $V_N$, namely if 
$$k(1+{2\over n}) = -\hbar(1 +2N/n),~~k'(1+{1\over n}) = \hbar^2{N\over
n}({N\over n} + 1).\eqno(5.10)
$$  
It turns out that such a
deformation exists if and only if the parameters $k$ and $k'$ are related to
each other,  precisely as implied by (5.10).

 \b
\def\p{\partial}

\no{\bf  5.4.3. Associated representation of
 so(n)   and sp(2n).}~
Let $V_N$ be as in 5.4.2, and $\eta$ a symmetric, nondegenerate 2-form on
$\Crm^n$. Let $\partial_a = \eta_{ab}\partial/\partial x_b$ and $\tilde L_{ab}
=
\hbar(x_a\partial_b - x_b\partial_a)$, which gives  a formal representation
of $so(n,\Crm)$ and the formula
$$
 \tilde L_{ab}\tilde L_{cd} =   \hbar^2(x_ax_c\p_b\p_d +\eta_{bc}x_a\p_d   -
a,b - c,d). 
$$
We want to simplify this representation as much
as possible and therefore restrict the variables to the cone $\eta(x,x) = 0$ 
and the space $V_N$ to the subspace of harmonic functions. (One  verifies 
that the choice $N = 2-n/2$ of the degree of
homogeneity   allows the Laplace operator a well defined action on 
functions defined on the cone.)
The first term on the right now satisfies the constraints.   
This is is not yet a model for a star product since the first order
operators on the right do not combine to
$\tilde L$'s. But another way to write the last relation is
$$\eqalign{
 \tilde L_{ab}\tilde L_{cd} = &~  \hbar^2\big(x_ax_c\p_b\p_d +
{1\over 2}(\eta_{bc}x_a\p_d + \eta_{ad}x_d\p_a) - a,b-c,d\big) \cr
& + {\hbar\over 2}\eta_{be}(L_{ad} - a,b-c,d).
 \cr}
$$
The first line satisfies the constraints and if we take this as a model for
the classical part then we are led to look for a star product of the form
$$
L_{ab}*L_{cd} = L_{ab}L_{cd} +   {\hbar\over 2}(\eta_{bc}L_{ad} - a,b-c,d) +
o(\hbar^2).
$$  
It turns out that such an invariant star product exist iff the last term is
precisely $\hbar^2(2-n/2)(1-n)$ times the Killing form. See [BZ] for a
complete discussion of singular representations of so$( p,q)$.

 In the case $n = 6, ~N = -1$, restricting to the real form
$so(4,2)$ and taking
$r_1,...,r_4\in \Nrm,~ r_5/r_6 \in \Nrm$, one recovers the unitarizable
representation that is realized on the space of solutions of a massless scalar
field in 4 dimensions, in the form discovered by Dirac. The same
representation appears in the theory of the hydrogen atom, where it is  
realized by self adjoint operators in $L_2(\Rrm,d^3x/r)$.

\b 
Let $V_{2n}$ be the space of polynomials on a complex vector space of
dimension $2n$ endowed with an antisymmetric, nondegenerate 2-form $\eta$ and
coordinates $\xi_1,...,\xi_{2n}$. A Poisson bracket is defined by
$\{\xi_a,\xi_b\} = \eta_{ab}$. Let $\partial_a:=\eta_{ab}\partial/\partial x_b$
and $\tilde L_{ab} = \hbar(\xi_a \partial_b + \xi_b\partial_a)$; then a very
similar analysis leads to the idea of a star product such that   
$$
 L_{ab}*L_{cd} = L_{ab}*L_{cd} + {\hbar\over 2}(
\eta_{bc}\tilde L_{ad}   + a,b + c,d) +  {\hbar^2\over 2}\eta_{ab}\eta_{cd}. 
$$
This is of course the Moyal star product, restricted to sp(2n).

\b\b\b

\ce{\steptwo 6.  Invariant star product on the minimal orbits.}

 \ce{\bf  6.1.  Computation of the homology.} ~A minimal, coadjoint
orbit of \g ~is an orbit through a highest weight vector of the coadjoint
action, equal to a highest root of \g ~. By a theorem
of Kostant [K2] one knows that the orbit is a variety
$$
M = \Cit[x_i,...,x_N]/(R),
$$
where $N$ is the order of \g ~ and (R) is the ideal generated by a set $R =
\{g_1,...,g_K\}$ of quadratic relations. We choose the quadratic
polynomials so that they are linearly independent 
and express  each $g_\alpha$ as $g_\alpha = g_\alpha^{ij}x_ix_j$ with
$g_\alpha^{ij}$ symmetric.
 
  We need a
generalization of the results in Section 3.2 to this case of multiple
relations,   up to the level of 3-chains and 3-cochains. It is clear
that Hoch$_2$ is the space spanned by the following chains,
$$
 (Z_{2,1})_{\alpha} = g_\alpha^{ij}x_i\otimes x_j,~~\alpha = 1,2,...~,~{\rm
and}~ ~~(Z_{2,2})_{ ij} = x_i\wedge x_j,~ i,j = 1,...,N~.
$$
Every closed 3-chain is homologous to a linear one,
$$
a = x_i\otimes x_j\otimes x_k A^{ijk}
$$
with $A^{ijk}\in \,\Cit$, and this chain is closed iff
$$
 x_ix_jA^{ijk} = 0 = x_jx_kA^{ijk} = 0.
$$
\vskip-2mm
\no Hence
$$
  A^{ijk} + A^{jik}  =    g_\alpha^{ij}c^k_\alpha,~~
 A^{ijk}+A^{ikj} =    c'^i_\alpha  g_\alpha^{jk},  
$$
with complex coefficients $c_\alpha,c'_\alpha$. This can be solved iff
$$
\sum_{cyc} g_\alpha^{ij}c^k_\alpha = \sum_{cyc} c'^i_\alpha  g_\alpha^{jk},
$$
where the sums are over cyclic permutations, which implies that  
  $c'^i_\alpha = c^i_\alpha + \rho^i_\alpha$, with the
coefficients $\rho_\alpha^i$ subject to $\sum_{cyc}
\rho_\alpha^ig_\alpha^{jk} = 0$,  as follows,
$$
6A^{ijk} = {\rm Alt}A^{ijk} +  3(g_\alpha^{ij}c^k_\alpha  
  - g_\alpha^{ik}c^j_\alpha + g_\alpha^{jk}c^i_\alpha)
+2(\rho^i_\alpha g_\alpha^{jk} - i,j),~~\sum_{cyc}
\rho_\alpha^ig_\alpha^{jk} = 0.
$$
The three terms lie in ${\cal Z}_{3,3}={\cal H}_{3,3},~{\cal Z}_{3,2}={\cal
H}_{3,2}$ and
${\cal Z}_{3,1}={\cal H}_{3,1}$,  respectively.  The first 
 
\no two are of the form
listed in Proposition 3.2.1. The third space, which is empty in the case that
there is only one relation, has not been determined.
\b

\no{\bf 6.1.1. Conjecture.}~{\it  The space ${\cal Z}_{3,1}={\cal
H}_{3,1}$ is spanned by   chains of the form} 
$$
e_3(1) g_\alpha^{ij}g_\beta^{kl}(x_i\otimes \{x_j,x_k\}\otimes x_l), 
$$
{\it where $\{g_\alpha\}$ is the full set of binary relations and $e_3(1)$ is
the BGS idempotent.}
\b
\no{\bf  6.1.2. Example.}~ For  sl($n$) the two relations
$U_a^bU_c^d - a,c = 0,~ U_d^fU_g^h - f,h = 0$ generate in this manner the
closed chain 
$$
U_a^b\otimes U_c^f\otimes U_g^h -a,c-f,h.
$$

The difficulty will be overcome with the help of the correspondence principle,
an adaptation of Weyl's symmetric ordering.

  \b
\ce{\bf  6.2. Invariant star products and correspondence
principle.}

We take a fresh point of departure. Suppose that an invariant star product
has the following property 
$$
S(x_{i_1}*...*x_{i_p}) = x_{i_1}...x_{i_p} +  \phi_{i_1...i_p},\eqno(6.1)
$$
where $S$ stands for symmetrization and the function $\phi_{\underline i}$ is a
formal series in
$\hbar$.
\b
\no{\bf  6.2.1. Remark.}~   We shall postulate that,  to each order
in
$\hbar$,
$\phi$ is a polynomial of order less than
$p$.   In this way we guarantee an
important property of the deformation: the Poincar\'e-Witt basis is preserved.
Actually, in the present context this is a weak limitation, since  
equivariant 2-cochains of higher order are scarce by reason of the constraints.
\b
Invariance of the star product imposes the requirement that the map
$\phi: A\rightarrow A$ defined by $\phi: x_{i_1}...x_{i_p}\mapsto
\phi_{i_1....i_p}$ be equivariant for the adjoint action. The only other
requirement is the obvious  one that the correspondence must be consistent
with the relations that define the variety.  Applying
these constraints to both sides of Eq.(6.1) results in conditions on the map
$\phi$. We shall calculate these conditions explicitly for monomials of order 
2 and 3. We shall show that these conditions are precisely the same as those
implied by associativity, confirming the rather obvious fact that associativity
is not a separate concern. So long as the correspondence is consistent with the
constraints (now the only issue), detailed   knowledge of the restricted
homology spaces is not required. Cohomology was crucial for the demonstration  
of an extension to higher orders (Theorem 2.1.3), but it is not the best tool
for establishing the basis at low orders.

Let the symbol $S$ stand for symmetrizationn in the order of factors and set
$$
S(x_i*x_j) = x_ix_j +  \psi_{ij},~~ S(x_i*x_j*x_k) = x_ix_jx_k + 
\phi_{ijk},\eqno(6.2)
$$
and recall that $x_i*x_j - x_j*x_i = \hbar\{x_i,x_j\} =\hbar
{\epsilon_{ij}}^mx_m$. The polynomials $\phi_{ijk}$ and $\psi_{ij}$ are assumed
to define equivariant maps as explained above.
\b 
\no{\bf 6.2.1. Proposition.}~ {\it Assume that the
equivariant polynomials $\psi$ and $\phi$ have been chosen so that relations
(6.2) are consistent with the constraints. Then there is an invariant,
associative star product such that (6.2) holds.}
\b
\no{\bf  Proof.}   ~ It follows from (6.2) that
$$
(x_ix_j)*x_k = x_i*x_j*x_k - {\hbar\over 2}{\epsilon_{ij}}^mx_m*x_k -
\psi_{ij}* x_k.
$$
The symmetrized star product is
$$\eqalign{
 S(x_i*x_j*x_k)&= {1\over 6}\sum_{\sigma\in S_3}(x_i*x_j*x_k)\cr &   = {1\over
2} x_i*x_j*x_k    + 
{ \hbar\over 3} x_i*\{x_k,x_j\}   +
{\hbar\over 6} \{x_k,x_i\}*x_j + i,j ,\cr}
$$
and we deduce that
$$\eqalign{ 
(x_ix_j)* x_k  =  &~x_ix_jx_k + \phi_{ijk}\cr & - \psi_{ij}* x_k - {\hbar\over
2}\Big(\psi(\{x_k,x_i\},x_j) + i,j\Big))\cr & - {\hbar\over
2}\Big( \{x_k,x_i\}x_j  + i,j\Big) - {\hbar^2\over
12}\Big(\{x_i,\{x_k,x_j\}\}  + i,j\Big).\cr} 
$$
Taking $\psi$ to be equivariant leads to some cancellations,
$$\eqalign{ 
(x_ix_j)* x_k  =  &~x_ix_jx_k + \phi_{ijk}  - \psi_{ij} x_k -
\psi(\psi_{ij},x_k) \cr & - {\hbar\over 2}\Big( \{x_k,x_i\}x_j  +
i,j\Big) - {\hbar^2\over 12}\Big(\{x_i,\{x_k,x_j\}\}  + i,j\Big).\cr}\eqno(6.3)
$$

Similarly 
$$\eqalign{ 
 x_i*(x_j x_k)  = &~x_ix_jx_k +  \phi_{ijk}  - x_i \psi_{jk}   -
-\psi(x_i,\psi_{jk}) \cr & - {\hbar\over
2}\Big( \{x_k,x_i\}x_j  + j,k\Big) - {\hbar^2\over
12}\Big(\{\{x_j,x_i\},x_k\}   + j,k\Big).\cr}
$$
These equations yield explicit expressions for the the values of the
two-forms $C_1$ and $C_2$  defined by 
$$\eqalign{
(x_ix_j)* x_k  = &  x_ix_jx_k + \hbar C_1(x_ix_j, x_k) + \hbar^2 C_2(x_ix_j,
x_k),\cr x_i*(x_jx_k)  =  &  x_ix_jx_k + \hbar C_1(x_i,x_jx_k) + \hbar^2
C_2(x_i,x_jx_k). 
\cr}\eqno(6.4)
$$ 
These values solve the condition for associativity of the star product on
linear chains. Theorem 2.3.1 then assures us that they can be satisfied in
general. The proposition is proved.

\b
In the present approach associativity
is satisfied trivially. What is far from trivial is the existence of a
function $\phi_{ijk}$ that solves (6.3). The obstructions are the
constraints. Application of
 $g_\alpha^{ij}$
to the first and $g_\alpha^{jk}$ to the second, gives
$$ 
g_\alpha^{ij}\phi_{ijk} = (g_\alpha^{ij}\psi_{ij})*x_k  + \hbar
g_\alpha^{ij}{\epsilon_{ki}}^m\psi_{mj} + {\hbar^2\over 6}
g_\alpha^{ij}{\epsilon_{ki}}^m{\epsilon_{mj}}^p x_p,\eqno(6.5)
$$
$$  
g_\alpha^{jk}\phi_{ijk} = x_i*(g_\alpha^{jk}\psi_{jk}) + \hbar
g_\alpha^{kj}{\epsilon_{ki}}^m\psi_{mj} + {\hbar^2\over 6}
g_\alpha^{jk}{\epsilon_{ki}}^m{\epsilon_{jm}}^p x_p . \eqno(6.6) 
$$
Since $\phi$ is symmetric both right hand expressions must agree,  
$$
(g_\alpha^{ij}\psi_{ij})*x_k - x_k*(g_\alpha^{ij}\psi_{ij}) = 
2\hbar g_\alpha^{ij}{\epsilon_{ik}}^m\psi_{mj}~.  
$$  
The only equivariant tensors available to use for the 2-chain $\psi$ are
the Killing form and, in the case of $sl(n)$, a term linear in the
generators, as above. Both satisfy this last condition, so of the two equations
(6.5-6) it is enough to examine the first. The problem of consistency of 
(6.1) is reduced to the existence of $\phi_{ijk}$ that solves 
Eq.s (6.5).

A complete determination of the restricted cohomology is not
required; it is enough to know the relations that define the orbit.

\b\b

\ce{\bf  6.3. Calculations for sl($n$).}

\no{\bf  6.3.1. Solving the constraints.}~Applied to sl(n), $x_i
\rightarrow  U_a^b$, with the notation and the commutation relations as in
Sections 5.2-4,
  the above result (6.3) take the form
$$
{1\over 2}(U_a^b*U_c^d + U_c^d*U_a^b) =  U_a^bU_c^d +  \psi_{ac}^{bd}
$$
and
$$\eqalign{
(U_a^bU_c^d)*&U_e^f =  U_a^bU_c^dU_e^f + \phi_{ace}^{bdf}  -
\psi_{ac}^{bd} U_e^f      - \psi(\psi_{ac}^{bd},U_e^f) 
 \cr &
\hskip1.5cm -{\hbar\over 2}\Big( \{U_e^f,U_a^b\} U_c^d  +
 \{U_e^f,U_c^d\} U_a^b \Big) \cr & \hskip2cm -{\hbar^2\over
12}\Big(\{U_a^b,\{U_e^f,U_c^d\}\} +
\{U_c^d,\{U_e^f,U_a^b\}\}\Big)  
\cr}\eqno(6.7)
$$
A better notation is to express $u\in \, $\g$\,'$ as $u = (AU) =
A^a_bU_a^b$. The oefficients
$A^a_b$ are coordinates for
\g$\,$; that is, $A$ ranges over the matrices of the adjoint
representation of \g. The first relation becomes
$$
u*u = u^2 + \psi(u,u).
$$
Equivariance restricts $\psi$:
$$
\psi(u,u) = k(AAU) + k'(AA),~~ k,k'\in \Cit,\eqno(6.8)
$$
where $(AAU)$ and $(AA)$ indicate traces of product of matrices. Similarly,
$$
u*u*u = u^3 + \phi(u,u,u),
$$
\vskip-.5cm
\no with
$$\eqalign{
\phi(u,u,u) = & \phi_1 (AAU)u + \phi_2(AAAU) +
\phi_3(AA)u +
\phi_4(AAA).
\cr} \eqno(6.9)
$$
(There are other invariants but their inclusion here is not allowed by the
relations.) Let $\psi$ be given in the form (6.8) and look at Eq.(6.7) as an
equation to determine $\phi$. Recall that the minimal orbit is defined by the
relations
$U_a^bU_c^d = U_a^dU_b^b,~ a,..,d = 1,...,n$.
\b
\no{\bf  6.3.2. Proposition.}~{\it Eq.(6.7) is consistent with the
relations $U_a^bU_c^d = U_a^dU_b^b$ if and only if the parameters $k$ and $k'$
are related to each other as follows,}
$$
4k'(1 + {1\over n}) = k^2(1+{2\over n})^2 - \hbar^2.\eqno(6.10)
$$
\b
The proof is in the Appendix.

\b\b
\def\d{\delta}
\no{ \bf 6.3.3.  Generators of the ideal. Proposition.}~{\it The
Joseph ideal for  sl($n$) is generated by the relations}
$$\eqalign{
U_a^b*U_c^d - U_a^d*U_c^b = &~(\hbar/2)(\d_c^bU_a^d - \d_a^dU_c^b
- \d_c^dU_a^b + \d_a^bU_c^d) \cr &
+ (k/2)(1+{2\over n})\big(\d_a^dU_c^b + \d_c^bU_a^d -  \d_a^bU_c^d -
\d_c^dU_a^b)\big)\cr &
+ k'(1 + {1\over n})(\d_a^d\d_c^b -  \d_a^b\d_c^d).
\cr}\eqno(6.11)
$$ 
\b
\no{\bf Proof.}~   
We find the relations of the deformed algebra by eliminating the original
product from 
$$\eqalign{
U_a^b*U_c^d = &~U_a^bU_c^d + (\hbar/2)(\d_c^bU_a^d - \d_a^dU_c^b) \cr &
+ (k/2)\Big(\d_c^bU_a^d + \d_a^dU_a^b - {2\over n}(\d_c^dU_a^b +
\d_a^bU_c^d)\Big)\cr &
+ k'(\d_a^d\d_c^b - {1\over n}\d_a^b\d_c^d).
\cr}\eqno(6.12)
$$
\b
We have seen that this relation is equivalent (given the Poisson bracket) to
the simple relation $u*u = u^2 + k(AAU) + k'(AA)$. 

A highest weight module over the deformed algebra  is a module generated by a
vector $v$ such that
$$
U_a^b*v = 0,~ a<b = 2,...,n,~~ U_a^a*v = \lambda_a v,~ a = 1,...,n,\eqno(6.13)
$$
with $\lambda \in  ~\Cit^n, ~\sum_1^n\lambda_a = 0$. (The Cartan subalgebra
consists of diagonal matrices.)
\b 
\no{\bf 6.3.4. Highest weight module. Proposition.}~{\it A highest
weight module of the deformed algebra exists if and only if the relation
(6.10) holds. In that case the highest weight is  one of the following, $1\leq
m\leq n$, 
$$
\lambda_1 = \lambda_2 = ... = \lambda_{m-1} =-{\hbar\over 2} - \gamma,~~ 
  \lambda_{m+1} = ... = \lambda_n = {\hbar\over 2} - \gamma .\eqno(6.14)
$$
with $\lambda_m$ determined by the fact that $\sum\lambda_i = 0$ and $\gamma = 
(k/2)(1+{2\over n})$.}
\b

\no{\bf Proof.}~The only relations that are changed by the
deformation are those where
$a,c$
 intersects $bd$, so we may limit ourselves to the case that $b = a$. Also, 
if $a$ is equal to $c$ or to $d$ the relations are just the usual
commutation relations; so take $c,d \neq a$. Then 
$$\eqalign{
U_a^a*U_c^d - U_a^d*U_c^a = &~{\hbar\over 2}(- \d_c^dU_a^a +  U_c^d)    
- {k\over 2}(1+{2\over n})\big(  U_c^d +
\d_c^dU_a^a)\big)    
- k'(1 + {1\over n})  \d_c^d.
\cr} 
$$
In particular, for $a\neq c$,
 $$
U_a^a*U_c^c - U_a^c*U_c^a =  ~(\hbar/2)(  U_c^c -  U_a^a)  
- (k/2)(1+{2\over n})\big( U_c^c  +U_a^a)\big) 
- k'(1 + {1\over n}).
$$
Applying this to the highest weight vector we obtain, for $a>c$,  
$$
\lambda_a\lambda_c = (\hbar/2)(\lambda_c-\lambda_a) -\gamma(\lambda_a +
\lambda_c) - \gamma',~~ \gamma := (k/2)(1+ {2\over n}),~ \gamma' :=
k'(1+{1\over n}).
$$
It follows that, 
$$
\lambda_a = {({\hbar\over 2} - \gamma)\lambda_c - \gamma'\over \lambda_c +
{\hbar\over 2} + \gamma},~~ a>c~.
$$
Suppose $n > 2$. Assume that $\gamma' \neq \gamma^2 - ({\hbar\over 2})^2$,
then taking $c = 1$ we get $\lambda_2 = ... = \lambda_n$ and taking $a = n$ 
we find that $\lambda_1 = ... = \lambda_{n-1}$.
Since $\sum \lambda_i = 0$ this is not interesting. We conclude that
 $$
\gamma' = \gamma^2 - ({\hbar\over 2})^2,\eqno(6.15)
$$
which is the same as (6.10).
Then there is an integer $m,~ 1\leq m \leq n$ such that $\lambda_m$
differs from its neighbours, and the statement of the proposition  follows
immediately.
\b

\no{\bf  6.3.5. Abelian deformation.}~ When $\hbar = 0$ the
deformed algebra is commutative, so a representation is just a character, or a
maximal ideal. The group acts on the maximal ideals and among the maximal
ideals   there are upper triangular ones. We have seen, Eq.(6.10) that this is
possible only if either $n = 2$ or if (6.10) holds  
$$
4k'(1 + {1\over n}) = k^2(1+{2\over n})^2.
$$
In this case the deformed variety is the space of traceless matrices with all
but one of the eigenvalues equal to $ -{k\over 2}(1+{2\over
n})$.   
\b 
 
\ce{\bf 6.4. Calculations for so(n).}
\no{\bf  6.4.1. Solving the constraints.}~Applied to sl(n), $x_i
\rightarrow  L_{ab}$, with the notation and the commutation relations as in
Section 5.4.3, we have 
$$
{1\over 2}(L_{ab}*L_{cd} + L_{cd}*L_{ab}) =  L_{ab}L_{cd} +  \psi_{ab,cd} 
$$
and   the above result (6.3) takes the form
$$\eqalign{
(L_{ab}L_{cd})*&L_{ef} =  (L_{ab}L_{cd})L_{ef} + \phi_{ab,cd,ef}  -
\psi_{ab,cd}  L_{ef}    
 \cr &
\hskip1.5cm -{\hbar\over 2}\Big( \{L_{ef},L_{ab}\} L_{cd}  +
 \{L_{ef},L_{cd}\}L_{ab} \Big) \cr & \hskip2cm -{\hbar^2\over
12}\Big(\{L_{ab},\{L_{ef},L_{cd}\}\} +
\{L_{cd},\{L_{ef},L_{ab}\}\}\Big)  
\cr}\eqno(6.16)
$$
A much better notation is to express $u\in  \,\g \,'$ as $u = (AL) =
A^{ab}L_{ab}$, the coefficients $A^{ab} = -A^{ba}$   coordinates for \g$\,$;
then these equations become 
$$
u*u = u^2 + \psi(u,u).\eqno(6.17)
$$
and
$$
u^2*v = u^2v + \phi(A,A,B) - \psi(A,A)v -  \hbar [v,u]u  
  -{\hbar^2\over 6}\{u,\{v,u\}\}. 
$$
 Equivariance restricts $\psi$ and $\phi$:
$$
\psi_{A,A} = k(AA),~~ 
 \phi(AAA) = \phi_1 (AAAL) + \phi_2(AA)u,  
$$
where $(AAU)$ and $(AA)$ indicate traces of product of matrices. (In this
context we take the form $\eta$ that defines  so($n$) to be
$\eta_{ab} = 0, \, a\neq b,\, \eta_{aa} = 1, a = 1,...,n$, so that
two-forms become matrices without any fuss.) 
We fix the parameter $k$ and look at Eq.(6.16) as an equation to
determine the coefficients $\phi_1,\phi_2$. Then we have:
\b
\no{\bf  6.4.2. Proposition.}~{\it Eq.(6.16) is consistent with the
relations that define the minimal orbit if and only if the parameter  $k$
takes the value $\hbar^2(n-4)/(n-1)$. In this case the parameters $\phi_1$ and
$\phi_2$ are fixed by (6.16).}
\b
The proof is in the Appendix.

\b

\no{\bf 6.4.3. Generators of the ideal. Proposition.}~{\it The
Joseph ideal for  so($n$) is generated by (the commutation relations
and) the relations}
$$
S(\eta^{bc}L_{ab}*L_{cd}) + {\hbar^2\over 2}(n-4)\eta_{ad}  = 0,~~
\sum_{cyc(abc)}(L_{ab}*L_{cd} -\hbar \eta_{ad}L_{bc})=
0.\eqno(6.18)
$$

These relations were derived  by Binegar and Zierau [BZ], who also
determined the highest weight module and   the associated unitary
representations of SO($p,q$). We are  interested, nevertheless, in deriving
these  results with the help of the star product. Eq.(6.18) is a direct
consequence of (6.17), that read in full
$$
L_{ab}*L_{cd} = L_{ab}L_{cd} + {\hbar\over 2}[L_{ab},L_{cd}] +
{k\over 2}(\eta_{ac}\delta_{bd} - a,b). 
$$ 
\b\b

\ce{\bf    6.5. Highest weight module for $so(2\ell+1)$.}~
Let
$\{E_{ab}\}, a,b = 1,..., $ denote the unit $n$-by-$n$ matrices as
earlier, and take
$$
\eta = \sum_{a = 1}^n\delta_{a,a'},~~ a\,' := n+1-a.
$$
A basis for so($\eta$) in the natural representation is the set $ \{L_{ab}\},~
a<b\in\{i,...,n\}$ of matrices
$$
L_{ab} = (E_{ab}- E_{ba})\eta = E_{ab'} - E_{ba'}.
$$
The Cartan subalgebra of choice has the basis
$$
H_a = L_{aa'}-L_{a'a},~~ a = 1,...,l.
$$
The set of positive root vectors is the collection
$\{L_{ab}\}, ~ a+b\leq n,~ a<b$.
Our calculations are insensitive to the parity of $n$; nevertheless we note
the following facts.  When $n = 2\ell + 1$  the simple roots are $L_{1,n-1},
L_{2,n-2},...,,L_{\ell-1,\ell},L_{\ell-1,\ell+1}$ and the associated roots are
$\alpha_i = (0,...,0,1,-1,0,...,0), ~i = 1,...,\ell-1$ with 1 in the $i$'th
place and
$\alpha_\ell = (0,...,0,1)$. When $n = 2\ell$ the simple roots are $L_{1,n-1},
L_{2,n-2},...,L_{\ell-1,\ell}, L_{\ell-1,\ell+1}$ and the associated roots are
$\alpha_i = (0,...,0,1,-1,0,...,0), ~i = 1,...,\ell-1$ with 1 in the $i$'th
place and
$\alpha_\ell = (0,...,0,1,1)$. All this as  in Bourbaki.

Having fixed the root system and the Cartan subalgebra we define a highest
weight module as one generated by a character $H_a \mapsto \lambda_a\in\, \Cit$
and a (highest weight) vector 
$v_0$ with the property
$$\eqalign{&
L_{ab}v_0 = 0, ~~ a + b \leq n,\cr &
H_av_0 = \hbar\lambda_av_0,~~   a = 1,...,\ell.
\cr}
$$

\b
\no{\bf Proposition 6.5.1.} ~ {\it A highest
weight module of the deformed algebra exists if and only if the parameter  $k$
takes the value $\hbar^2(n-4)/(n-1)$. In that case the highest weight is  one 
of the following,   for}   $k\in \{1,...,\ell\}$,
$$
\lambda_1 = ...= \lambda_{k-1} = -1,~~ \lambda_k = k + 1 -{n\over 2}, ~~
 \lambda_{k+1} = ... = \lambda_\ell = 0.
$$
\b
The proof  is in the Appendix.

  Joseph's choice is $ k = \ell-1$. When $n = 2\ell+1$, so($n$) = so(2$\ell+1$)
=
$B_\ell$, 
$$\eqalign{
\lambda + \rho &= \omega_1 + ... + \omega_{\ell-3} + {1\over
2}(\omega_{\ell-2} + \omega_{\ell-1}) + \omega_\ell\cr & = (\ell -
{3\over 2},\ell-{1\over 2},...
 ,{5\over 2},{3\over 2},1,{1\over 2}),~~ \lambda = (-1,...,-1,-{1\over 2},0).
\cr}
$$
When $n = 5$ this corresponds to the Bose singleton. Another choice is
$\lambda = (-1,..,-1,{1\over 2})$; in the case $n = 5$ it is the Fermi
singleton. Both are unitarizable (after taking a quotient)
representations of so(3,2). Binegar and Zierau take
$k = 1$.

When $n = 2\ell,$ so($n$) = so(2$\ell$) = $D_\ell$, Joseph again takes $k =
\ell -1$ and
$$\eqalign{
\lambda + \rho &= \omega_1 + ... + \omega_{\ell-3} + \omega_{\ell-1} +
\omega_\ell\cr 
&
= (\ell-3,\ell-2,...,1,1,0),~~ \lambda = (-1,-1,...-1,0,0).
\cr}
$$
When $n = 6$ this is the highest weight  
of the representation of the conformal group by a scalar massless field in 4
dimensions, the same representation that appears in Schroedinger's hydrogen
atom. Binegar and Zierau take $k = 1$.
      
\b 

\ce{\bf  6.6. Calculations for  the exceptional simple Lie
algebras.}

\no{\bf  6.6.1. The relations.}~ Let $K$ denote the Killing form,
$$
K(x_i,x_j) = K_{ij} = -{\rm tr}(x_ix_j) = -{\epsilon_{mi}}^n{\epsilon_{nj}}^m
$$ and $K^{ij}$ the matrix elements of the inverse matrix.
 The reduction of the
symmetric part of the adjoint representation is governed by the operator $L$in
\g~$\otimes$ \g ~ defined by 
$$
L: x_i\otimes x_j \mapsto  K^{mn}{\epsilon_{im}}^s{\epsilon_{jn}}^tx_s\otimes
x_t = K^{mn}[x_i,x_m]\otimes [x_j,x_n].\eqno(6.19)
$$
The symmetric product of the adjoint representation by itself decomposes into
a sum of three terms, $Id \oplus V_2 \oplus V_3$. The last contains the
extremal weights and is the minimal orbit in \g. The relations that define the
minimal orbit are the projections on the first two. The one-dimensional
component expresses the condition that 
$$
K^{ij}x_ix_j = 0,
$$
jointly they fix  the eigenspace of $L$ in \g~$\otimes\,$\g. 
The  eigenvalues $l_i = l_1,l_2,l_3 $ associated with the eigenspaces
$Id, V_2,V_3$ are as follows. In all cases, $l_1$ = 1, and $l_2 + l_3 =
-1/6$. The values of $l_3$ are taken from [MP],
$$
\matrix{\g&G_2&F_4&E_6&E_7&E_8\cr
l_3&1/4&1/9&1/12&1/18&1/30\cr D & 27&52&78&133&248} 
$$
The number $D$ is the dimension of \g.

We turn to Eq.(6.3). Actually, the exceptional simple Lie algebras are the
easiest to deal with, and all five can be done uniformly. The immediate reason
for this is the non-existence of an irreducible, invariant fourth order
polynomial. The  only
equivariant, symmetric tensors are (see 6.2.1)   
$$
\psi_{ij}  = kK_{ij},~~ \phi_{ijk} = {k'\over3}\sum_{cyc}K_{ij}x_k,\eqno(6.20)
$$
with $k$ and $k'$ to be determined by the relations. So that Eq.(6.3)
simplifies: 
  
$$\eqalign{ 
(x_ix_j)* x_k & -  x_ix_jx_k + {\hbar\over 2}\Big( \{x_k,x_i\}x_j  +
i,j\Big)\cr &=  \phi_{ijk}  - \psi_{ij} x_k -
     {\hbar^2\over 12}\Big(\{x_i,\{x_k,x_j\}\}  +
i,j\Big).\cr}\eqno(6.21) 
$$
 
\b
\no{\bf  6.6.2. Theorem.}~{\it Let  \g~ be one of the five
exceptional simple Lie algebras, and $K$ the Killing form. There exists a 
unique, invariant star product on the minimal nilpotent orbit, such that
$$
S(x_i*x_j) = x_ix_j + kK_{ij},~~ k \in \Cit,
$$
for one and only one choice of $k, k'\in \, \Cit$, namely $k   
  = {l_2-1\over D}{\hbar^2\over 6} + l_2{\hbar^2\over 12},~ k' =
{\hbar^2\over 4}$.}
\b

\no{\bf  Proof of the Theorem.}~  
  It has already been pointed out that  
the only possible forms of the deformed products $x_i*x_j = x_ix_j + o(\hbar)$
and $S(x_i*x_j*x_k) = x_ix_jx_k + o(\hbar)$ are
$$
S(x_i*x_j) = x_ix_j + \psi_{ij}x_k,~~ S(x_i*x_j*x_k) = x_ix_j + \phi_{ijk},
$$
with the equivariant maps $\psi,\phi$ as in (6.20). 
The left side of (6.21) satisfies the constraints; the
question is whether the right side does. Fix the \underbar{index} $_k$ and
define vectors
$X,Y$ in
\g~$\otimes\,
$\g~ with components
$$
X_{ij} =  K_{ik}x_j + i,j,~~ Y_{ij} = K_{ij} x_k.
$$
    Let $L_s$ denote the projection of the operator $L$ on
the symmetric part of 
\g~$\otimes\,$\g$\,$; then the right hand side of (6.21) is the vector
$$
{k'\over 3}( X+Y) - kY - {\hbar^2\over 12}L_sX.\eqno(6.22)
$$ 
The operator $L_s$ satisfies the characteristic equation
$(L_s-1)(L_s-l_2)(L_s-l_3) = 1$, and since $L_sY = Y$ there is a constant $c$
such that 
$$
(L_s-l_2)(L_s-l_3)X = cY.
$$
To determine the constant we contract this equation with $K$ and find
that
$$
c = {2\over D} (1-l_2)(1-l_3).
$$
With this, the vector (6.22) reduces to
$$
Z := c^{-1}({k'\over 3}-k)(L_s-l_2)(L_s-l_3)X + {k'\over 3}  X -{\hbar^2\over
12}L_sX
$$ 
Projecting on the trivial representation we get
$$
0 =  {k'\over 3} (D + 2) -kD  - \hbar^2/6 = 0,~~ D = {\rm dim}~
\g\,.\eqno(6.23)
$$
All the constraints are expressed by $(L_s-l_3)Z = 0$, which yields 3
conditions
$$\eqalign{
(1-l_3)c^{-1}({k'\over 3} -k) -{\hbar^2\over 12} &= 0,\cr
 (1-l_3)c^{-1}({k'\over 3} -k)(l_2 + l_3) + {k'\over 3} &=
   l_3{\hbar^2\over 12},\cr
(1-l_3)c^{-1}({k'\over 3} -k)l_2l_3 &=  
l_3{k'\over 3}.
\cr}\eqno(6.24)
$$ 
All four equations (6.23-4) agree on unique values of $k$ and $k'$,
$$
k' = l_2{\hbar^2\over 4},~~ k -{k'\over 3} = {l_2-1\over D}{\hbar^2\over 6}.
$$
The proposition is proved.

\b
\no{\bf  6.6.3. The Josef ideal.}~ It is of course generated by
the relations
$$
x_i*x_j - x_j*x_i = \hbar\{x_i,x_j\},
$$
and
$$
\Big(L_s-l_3\Big))(x_i*x_j-
kK_{ij}) = 0.
$$ 

\b\b

\no{\steptwo Acknowledgements }

I thank Birne Binegar for information about the Joseph ideals.

\ve

\no{\steptwo  Appendix}

\no{\bf Proof of Proposition 6.3.2.}~ The lack of symmetry in (6.7)
forces the introduction of a second variable,
$v = (BU)$; then the equation can be written as
$$
u^2*v = u^2v+\phi(u,u,v) -\psi(u,u)v -  \psi(\psi(u,u),v) 
-\hbar \{v,u\}u  - {\hbar^2\over 6} \{u,\{v,u\}\}, 
$$ 
with $\phi$ as in (6.9),
$$\eqalign{
\phi(u,u,v) = &  
+ {\phi_1\over 3} \Big((AAU)v + 2(ABU)u\Big)      
+ {\phi_2\over 3}\Big((AABU) +(ABAU) + (BAAU)\Big)\cr & +
{\phi_3\over 3}\Big((AA)v + 2(AB)u\Big) + 
\phi_4(AAB).
\cr} 
$$
Exchange of the upper indices on the two $A$'s is a type of Fierz
transformation. The effect is as follows.  
$$\eqalign{
  (ABAU) \mapsto & ~(AB)u - {1\over n}(AABU)- {1\over n}(AAUB)+ {1\over
n^2}(AA)v\cr
(AABU)\mapsto & -{2\over n}(AABU) + {1\over n^2}(AA)v \cr 
  (ABU)u \mapsto &~(ABU)u - {1\over n}(AAU)v\cr
(BAU)u \mapsto &~ (BAU)u - {1\over n} (AAU)v,\cr
(AAU) \mapsto & -{2\over n}(AAU),~~~~~(AAB) \mapsto  -{2\over n}(AAB),\cr
(AB)u \mapsto &~(AUAB),~~~~~
~~~~~\,(AA) \mapsto  -{1\over n}(AA). 
\cr}
$$
and 
$$\eqalign{
\phi(u,u,v) \mapsto &  
+ {\phi_1\over 3} \Big(-{2\over n}(AAU)v + 2(ABU)u -{2\over n}(AAU)v\Big)     
\cr & + {\phi_2\over 3}\Big(-{3\over n}(AABU + AAUB) +{3\over n^2}(AA)v  
   +(AB)u \Big)\cr
& +{\phi_3\over 3}\big(-{1\over n}(AA)v + 2(ABAU)\big)  
 +  \phi_4{-2\over n}(AAB). 
\cr} 
$$
Furthermore,
$$\eqalign{ 
-\psi(u,u)v &= -k(AAU)v  - k'(AA)v\mapsto {2k\over n}(AAU)v + {k'\over
n}(AA)v\cr   -  \psi(\psi(u,u),v) &=   -{k^2\over 2}\Big((AABU + AAUB)
-{2\over n}(AA)v\Big)  - kk'(AAB)\cr &\hskip1cm\mapsto   {k^2\over
n}\Big((AABU + AAUB) -{2\over n}(AA)v\Big)\cr &\hskip1in -{k^2\over n}(AA)v +
{kk'\over n}(AAB)\cr     - {\hbar^2\over 6}
\{u,\{v,u\}\} &=- {\hbar^2\over 3} (ABAU) + {\hbar^2\over 6}(AABU + AAUB)\cr  
\hskip1in &\hskip1cm \mapsto  {\hbar^2\over 3}\Big(- (AB)u + {1\over n}(AABU +
AAUB) -{1\over n^2}(AA)v\Big)\cr &
\hskip1in-{\hbar^2\over 6}\Big(-{2\over n}(AABU + AAUB) + {2\over
n^2}(AA)v.\Big)
\cr}
$$
 Applying the constraint $U_a^bU_c^d = U_a^dU_c^b$ to (6.7) we get      
$$\eqalign{
{1\over 3}\phi_1(1+{2\over n}) - k(1+{2\over
n}) = &~0,\cr  
2\phi_2(1+{3\over n}) -3k^2(1+{2\over n}) +  \hbar^2  =  & ~0,\cr   
 \phi_2 -2\phi_3-\hbar^2  =&~ 0,\cr  
\phi_3(1+{1\over n})-3k'(1+{1\over n}) -{3\over n^2}\phi_2 + {3\over
n}k^2(1+{2\over n}) = &~0,\cr
\phi_4(1+{2\over n}) -kk'(1+{2\over n}) = &~0.    
\cr}
$$
Eliminating the parameters $\phi_i$, we find the relation
$$
4k'(1+{1\over n}) = k^2(1+{2\over n})^2 - \hbar^2. 
$$
The proposition is proved.

\b\b

\no{\bf Proof of Proposition 6.4.2.}~  Recall that the relations are
$$
\sum_{cyc(abc)}L_{ab}L_{cd} = 0,~~ \eta^{bc}L_{ab}L_{cd} = 0. 
$$
Applying the first relation, 
 contraction on $b,c$, to Eq. (6.16) one gets
$$
0 = \eta^{bc}\phi_{ab,cd,ef} - {k\over 2}(n-1)\eta_{ad}L_{ef} 
-{\hbar^2\over 12}\Big((n-4)(\eta_{fd}L_{ae} + a,d - e,f)
-4\eta_{ad}L_{ef}\Big),
$$ 
and 
$$  
\eta^{bc}\phi_{ab,cd,ef} =    \Big({\phi_1\over 24}(n-1)- {\phi_2\over
6}\Big)(\eta_{fd}L_{ae} + a,d - e,f) 
  + \Big(-{\phi_1\over 6} + {\phi_2\over6}(n-1)\Big)\eta_{ad}L_{ef}.
$$ 
The second relation (a kind of Fierz transformation) gives
$$
0 = \sum_{cyc}\phi_{ab,cd,ef} + {\hbar^2\over 3}\eta_{fa}(\eta_{de}L_{bc} +
\eta_{ec}L_{db})-e,f 
$$
and
$$\eqalign{&
\sum_{cyc}\phi_{ab,cd,ef} = \big({\phi_1\over 12} +
{\phi_2\over 6}\big)\sum_{cyc}\eta_{fa}(\eta_{de}L_{cb} +
\eta_{ce}L_{bd})-e,f. 
\cr}
$$
Thus 
$$\eqalign{   
{\phi_1\over 24}(n-1)- {\phi_2\over 6} -{\hbar^2\over 12}(n-4)=& ~0,\cr
 -{\phi_1\over 6} + {\phi_2\over6}(n-1) - {k\over 2}(n-1) + {\hbar^2\over 3} =
&~0,\cr
 ~{\phi_1\over 12} +
{\phi_2\over 6} + {\hbar^2\over 3} = &~0.
\cr}
$$
Finally,
$$\eqalign{
\phi_1(n+1) = &~ 2\hbar^2(n-8),~~ \phi_2(1+n) = -3\hbar^2(n-2),\cr  k(n-1) =&~
\hbar^2(n-4).
\cr} 
$$
The proposition is proved.
\b\b

\no{\bf Proof of Proposition 6.5.1.}~All Cartan subalgebras are
isomorphic and any two systems of simple roots are related by a transformation
of the Weyl group.   The ideal determines only the infinitesimal character
$\chi(\lambda) := \rho +
\lambda$, where
$\rho$ is  half the sum of the positive roots and $\lambda$ is the highest
weight, up to a Weyl reflection. It does not distinguish between weights that
are related by a Weyl reflection of
the infinitesimal character.
  In the
case of $B_l = ~$so$(2\ell + 1)$ the formula is $\rho = (\ell -{1\over 2}, l-
{3\over 2},...,{1\over 2})$.

The   problem is to determine the possible values of the
infinitesimaal character.

We begin with the relation
$$
\sum_bL_{db}*L_{b'a} - {\hbar\over 2}(n-2)L_{da} + {k\over 2}(n-1)\eta_{da} =
0,
$$
in the case $a+d= n+1, a<d$. Applied to $v_0$ it gives $(L_{da} \rightarrow -
\lambda_a$)
$$
\Big(\sum_{b\leq a}L_{db}*L_{b'a} + {\hbar^2\over 2}(n-2)\lambda_a + {k\over
2}(n-1)\Big)v_0 = 0.
$$
Why? Well if $b>a$  then $b'<a'$ and $a + b' < n+1$ and $L_{b'a}v_0 = 0$.
Thus
$$
\Big(\sum_{b< a}\{L_{db},L_{b'a}\} + \lambda_a^2 + {\hbar^2\over
2}(n-2)\lambda_a + {k\over 2}(n-1)\Big)v_0 = 0,
$$
The bracket is $\hbar^2(\lambda_b - \lambda_a),~ b  = 1,...,a-1,$, so finally
$$
 (1-a)\lambda_a +  (\lambda_1 + ... + \lambda_{a-1}) + \lambda_a^2
+ {1\over 2}(n-2)\lambda_a + {1\over 2}(n-4)  = 0,
$$
In particular
$$
 \lambda_1^2 +{1\over 2}\lambda_1(n-2) + {1\over 2}(n-4) = 0,
$$
or
$$
(\lambda_1+1)\big(\lambda_1 + {n-4\over 2} \big) = 0.
$$
Thus
$$
  (1-a)\lambda_a +  (\lambda_1 + ... + \lambda_{a-1}) + \lambda_a^2
-\lambda_1^2 + {n-2\over  2} (\lambda_a-\lambda_1) = 0,
$$
and
$$
(\lambda_{a} -\lambda_{a-1})\Big(\lambda_a + \lambda_{a-1} + 
{1\over 2}(n-a)\Big) = 0.
$$

 We return to
$$
\sum_bL_{db}*L_{b'a} - {\hbar\over 2}(n-2)L_{da} + {k\over 2}(n-1)\eta_{da} =
0,
$$
now in the case that $a + d = n+2,~ 2\leq a<d$ . Applying to the highest
weight 
$$
\Big(\sum_{b\leq a}L_{db}*L_{b'a} - {\hbar\over 2}(n-2)L_{da}  
 \Big)v_0  = 0,
$$
or
$$
\Big(\sum_{b\leq a-2 }L_{db}*L_{b'a} + L_{da}*L_{a'a}  + L_{dd'}L_{da} -
{\hbar\over 2}(n-2)L_{da}  
 \Big)v_0  = 0.
$$
The first term can be replaced by the bracket $\sum_{b\leq a-2
}[L_{db},L_{b'a}] = \hbar\sum_{b\leq a-2 }L_{da}$, so   for $a =
2,...,\ell$,
$$
\Big( (a-1)\L_{da} -\lambda_aL_{da} - \lambda_{a-1}L_{da}-
{1\over 2}(n-2)L_{da} \Big)v_0  = 0.
$$
Hence
$$
\Big(\lambda_a + \lambda_{a-1} + {1\over 2}(n-2a)\Big)L_{da}v_0 = 0, ~~a =
2,...,
\ell.
$$
Note that $\{L_{da},L_{a'd'}\}v_0  =  \hbar^2(\lambda_a - \lambda_{a-1})v_0$.
The information contained in this last result is therefore precisely the same
as in (6.*), for when $\lambda_2=...=\lambda_n = 0$, then $L_{ab}v_0 = 0, ~a,b
> 1
 $.

Next, the other relation,
$$
\sum_{cyc bcd}(L_{ab}-\hbar \eta_{ab})*L_{cd} = 0,
$$
applied to $v_0$ 
in the case $a + b = c+d = n+1,~ a<b,c<d,a\neq c$, 
$$
 \Big((\lambda_a-1)\lambda_c +  L_{ac}*L_{db} + L_{ad}*L_{bc}\Big)v_0  = 0.
$$
Evaluating this in two cases, $a >b, a<c$ we find that $\lambda_<$  is -1 or
$\lambda_>$  is 0. That completes the proof.

\b\b

\no{\steptwo References.}~

\no [AA] Agarwal, G.S. and Wolf, E., Phys. Rev. {\bf D2} (1970) 2161. 

\no [AM] Abellanas, L. and Martinez-Alonzo, L., J. Math. Phys. (1976) 1363. 

\no [B]  Barr, M.,  Cohomology of commutative algebras,  
Dissertation, U.Penn. 1962;

    Harrison homology, Hochschild homology
and triples, J.Algebra {\bf 8} (1968)
 314-323.

\no [BF] Bayen, F., and Fronsdal, C., Quantization on the sphere, 

J. Math.
Phys. {\bf 22} (1981) 1345-1349.

\no [BFLS] F. Bayen, M. Flato, C. Fronsdal, A. Lichnerowicz and D. Sternheimer,

  Quantum  Mechanics as a Deformation of  Classical
Mechanics,

   Annals of Physics {\bf 111} (1978 ) 61- 110 and
111-151.

\no [BGS] Beilinson, A., Ginsburg, V. and Schechtman, V., Koszul Duality, 

J.
Geom. Phys. {\bf 5} (1988) 317-350.

\no [Bn] Berezin, F.A. General concept of quantization,

 Commun. Math. Phys.
{\bf 40} (1975) 153-174.
 \ve
\no [Be] Bezrukavnikov, R., Koszul property and Frobenius splitting of Schubert
varieties, 

 q-alg preprint, 9502021.

\no [BZ] Binegar, B. and Zierau, R., Unitarization of a singular
representation of SO($p,q$), 

Commun. Math. Phys. {\bf 138} (1991) 245-258.

\no [BJ] Braverman, A. and Joseph, A., The minimal realization from Deformation
Theory,

J. Algebra, {\bf 205}, 13-16 (1998).

\no [dL] de Wilde, M. and Lecomte, P.B.A., Existence of star products and of
formal 

deformations of the Poisson algebra of arbitrary symplectic manifolds,

Lett. Math. Phys. {\bf 7} (1983) 487-496.

\no [F]  Fedosov, B.V.,   A simple geometrical construction of
deformation quantization,

J.Diff.Geom. {\bf 40} (1994) 213-238.

\no [FG] Fronsdal, C. and Galindo, A., The ideals of free differential
algebras. 

J. Algebra, {\bf 222} (1999) 708-746.  

\no [FK] C. Fronsdal and M. Kontsevich, Quantization on Curves,
Math-ph/0507021. 
 
\no [FLS] M. Flato, A. Lichnerowicz and D. Sternheimer, Deformations of
Poisson brackets, 

Dirac brackets and applications, J. Math. Phys., {\bf 17}
(1976) 1754-1762.  

\no [Fy] Fleury  P.J., Splittings of Hochschild's complex for commutative
algebras,

Proc.AMS, {\bf 30} (1971) 405-323.

\no [F1] Fronsdal C., Some ideas about quantization, Rep. Math. Phys.
{\bf 15} (1978) 113.

\no [F2] Fronsdal C.,   Abelian Deformations, Proceedings of the IX'th
International  

ference on Symmetry Methods in Physics, Yerevan, July 2001.

\no [G] Garsia A.M., Combinatorics of the free Lie algebra and the
symmetric group,

Analysis, Et Cetera, Research Papers Published in Honor
of Jurgen Moser's 60'th

Birthday, Academic Press, New York, (1900)
309-382.

M. Barr's 60'th birthday, June 30, 1998.

\no [GS] Gerstenhaber M. and Schack S.D., A Hodge-type decomposition for
commutative

algebra cohomology, J. Pure and Applied Algebra
{\bf 48} (1987) 229-247.

\no [Gt]   Gutt, S., Lett. Math. Phys. {\bf 7} (1983) 249-258.

\no [G1] Gerstenhaber M., On the deformations of rings and algebras,

Annals of Math.{\bf 79} (1964) 59-103.
  
\no [G2] Gerstenhaber M., Developments from Barr's thesis, presented at
the celebration of 

the  60'th birthday of M. Barr.

\no [H] Harrison D.K., Commutative algebras and cohomology,

Trans.Am.Math.Soc. {\bf 104} (1962) 191-204.

\no [HKR] Hochschild G., , Kostant B. and Rosenberg A., Differential
forms on regular

affine algebras, Trans. Amer. Math. Soc. {\bf 102}
(1962) 383-408.

\no [J] Joseph, A., The minimal orbit in a simple Lie algebra and its
associated maximal 

ideal, Ann. Sci. Ecole Norm. Sup. {\bf 9},1-30 (1976).

\no [Kh ] Kontsevich M., Deformation quantization of Poisson
manifolds,
q-alg/9709040;

Operads and motives in deformation
quantization,

   Lett.Math.Phys. {\bf 48} (1999) 35-72.

  \no [K1] Kostant B., {\it  Quantization and unitary representations,
 Lectures in Modern

 Analysis and Applications III}, Lecture Notes in
 Mathematics {\bf 170},   87-208,

    Springer-Verlag, Berlin, 1970.

\no [K2] Kostant B., reported by A. Joseph.

\no [L] Loday,  J.-L. {\it Cyclic homology}, Springer-Verlag Berlin
Heidelberg, 1998.

\no [Ll]   Lledo, M.A., Deformation Quantization of Non Regular Orbits of
Compact Lie 

Groups, Lett.Math.Phys. 58 (2001) 57-67.

\no [M] Moyal, J.E.,  Quantum mechanics as a statistical theory,

 Proc.Cambridge Phil. Soc.{\bf 45} (1949) 99-124.

 \no [S] Souriau J.M., {\it Structures des Syst\`emes Dynamiques},
  Dunod, Paris, 1970.

\no [T] Tamarkin D.E., Another proof of M. Kontsevich' formality
theorem for $\Rrm^n$,

    math.QA/9803025.
 
\no [V] Vey J.,  D\'eformation du crochet de Poisson sur une vari\'et\'e
symplectique,

Comment. Math. Helv.  {\bf 50} (1975) 421--454.

\no [W]  Weinstein, A., Seminaire Bourbaki n$^o$ 789, 1993-94.

\no [Wl] Weyl H., {\it  Theory of Groups and Quantum Mechanics},
 Dover,
 New York, 1931.

\no [Wr] Wigner, E.P., Quantum corrections for thermodynamic equilibrium, 

Phys.
Rev. {\bf 40} (1932) 749.

\no [WS] Wee Teck Gan and G. Savin, Uniqueness of the Joseph ideal, 

found on
the home page of either author, U. Utah. Publication unknown.
\ve

\end

\no{\bf  'Proof' of conjecture.}~
Let $V$ denote the enveloping algebra of \g~ as a vector space, and
$V\rightarrow SV$ the canonical bijection on the symmetric algebra. We are
interested in subspaces
 $W_{12},~ W_{23}$ of $W := V^{\otimes 3}$ with the property that the mapping
$V^{\otimes 3}\rightarrow (SV)^{\otimes 3}$, followed by either of the two
the multiplication maps
$M_{12}: a\otimes b\otimes c\mapsto ab\otimes c,~ M_{23}: a\otimes b\otimes c
\mapsto a\otimes bc$ leads to zero.  More precisely, we need to find the
intersection of these two subspaces.

We have already found two thirds of this intersection, $Z_{3,3}$ arises from
commutativity, from the fact that $M_{12}$ annihilates $a\wedge b\otimes
c$, without any assistance from the ideal $R = (\{g_\gamma\})$.
The  chains in  ${\cal Z}_{3,2}$ are sums of terms that are killed by
multiplication, either by commutativity or by the vanishing of the $g_\gamma$,
these chains are linear in the $g_\gamma$. What is missing is the part that
depends on the vanishing of the $g_\gamma$ alone. It is clear from the
constraint 
$\sum_{cyc}\rho_\gamma g_\gamma^{ij} = 0$ that the $\rho$ are not independent
of the $g_\gamma$. This
heuristic argument suggested the following.

Consider the following operators in $W = V^{\otimes 3}$:
$$
(Q_{12})_\gamma   = g_\gamma^{ij}(x_i\otimes x_j\otimes 1)
$$
Each of these operators projects $W$ into the subspace $W_{12}$ and the union
of their images spans the part of this subspace that is symmetric in the two
first factors. Similarly define
$$
(Q_{23})_\gamma   = g_\gamma^{kl}(1 \otimes x_k\otimes x_l),
$$
with similar properties relative to the subspace $W_{23}$. We are interested
in the intersection of these two subspaces, and here we are fortunate.
\b
\no{\bf  6.4.4. Lemma.}~ {\it The intersection of the two
spaces is precisely the span of the images of the operators}
$$
 (Q_{12})_\gamma (Q_{23})_\beta  = g_\gamma^{ij}g_\beta^{kl}(x_i\otimes
 x_jx_k \otimes x_l) 
$$
and
$$
 (Q_{23})_\beta (Q_{12})_\gamma   = g_\gamma^{ij}g_\beta^{kl}(x_i\otimes
 x_kx_j \otimes x_l). 
$$

\b
\no The proof of the lemma consists of the observation that the
operator
$$
 [(Q_{12})_\gamma, (Q_{23})_\beta]  = g_\gamma^{ij}g_\beta^{kl}(x_i\otimes
 [x_j,x_k] \otimes x_l).
$$
projects into the intersection of the two subspaces; that is, the
commutator projects on a subspace of closed chains. Thus multiplication on the
first 2 factors gives
$$
2g_\gamma^{ij}g_\beta^{kl}(x_i[x_j,x_k]) = g_\beta^{kl}\{g_\gamma, x_k\} = 0.
$$
That completes the proof of the Lemma.  Only the
commutator is of interest, since it is linear and not exact. The proof of 
Theorem 6.4.3 is complete.
 \b
\no{\bf  6.4.4. Remark.} ~In the case that there is only one
constraint  the closed 3-chain constructed here is alternating and is
annihilated by the BGS idempotent $e_3(1)$. 
\b